\documentclass[11pt,a4paper]{article}
\usepackage{amsmath,amsthm,amsfonts,amssymb,mathrsfs}
\usepackage{graphicx}
\usepackage{subfigure}
\usepackage[latin1]{inputenc}
\usepackage[swedish,english]{babel}
\renewcommand{\vec}[1]{
{\boldsymbol#1}}

\begin{document}
\title{An explicit kernel-split panel-based Nyström scheme for
    integral equations on axially symmetric surfaces}
\author{Johan Helsing\thanks{Centre for Mathematical Sciences, Lund
    University, Sweden}~~and Anders Karlsson\thanks{Electrical and
    Information Technology, Lund University, Sweden}}
\date{\today}
\maketitle

\begin{abstract}
  A high-order accurate, explicit kernel-split, panel-based,
  Fourier--Nyström discretization scheme is developed for integral
  equations associated with the Helm\-holtz equation in axially
  symmetric domains. Extensive incorporation of analytic information
  about singular integral kernels and on-the-fly computation of nearly
  singular quadrature rules allow for very high achievable accuracy,
  also in the evaluation of fields close to the boundary of the
  computational domain.
\end{abstract}

\section{Introduction}

This work is on a high-order accurate panel-based Fourier--Nyström
discretization scheme for integral equations associated with the
Helmholtz equation in domains bounded by axially symmetric surfaces.
Efficient axisymmetric solvers for wave propagation and scattering are
important in their own right in optical and microwave
applications~\cite{Keka11,Wang08}. They are also needed in
multi-particle contexts, for example, to predict the effects of
absorption and scattering of sun light from soot in the
atmosphere~\cite{Nosi12}.

The present scheme resembles that of Young, Hao, and
Martinsson~\cite{Youn12}. The main difference lies in the treatment of
nearly- and weakly singular oscillatory kernels. In~\cite{Youn12}, a
precomputed 10th order accurate general-purpose Kolm--Rokhlin
quadrature~\cite{Kolm01} is used for discretization in the polar
direction and the post-processor, where field evaluations are done,
does not address the nearly singular case. Here we use instead 16th
order analytic product integration and a splitting of transformed
kernels. Quadrature weights are computed on the fly whenever needed.
While it has been considered hard to implement such an ``explicit
split'' panel-based scheme for axisymmetric problems,
see~\cite{Hao14}, our work demonstrates that it is indeed worth the
effort.

We specialize to the interior Neumann problem and are particularly
interested in finding solutions corresponding to homogeneous boundary
conditions (Neumann Laplace eigenfunctions). This PDE eigenvalue
problem models acoustic resonances in sound-hard voids and is of
importance in areas such as noise reduction~\cite{Huang02}, resonance
scattering theory~\cite{Colt13,Uberall82}, and quantum
chaos~\cite{Backer02,Barn06b,Barn14}. A related, vector valued,
problem models axially symmetric electromagnetic
scattering~\cite{Gedney90}. More recent general developments on
electromagnetic integral equation formulations can be found
in~\cite{Vico13}.

Aside from improving the convergence rate, our scheme improves the
achievable accuracy to the point where it can be called nearly
optimal. Furthermore, high accuracy is not only obtained for the
solutions to the integral equations under consideration. The
flexibility offered by performing weight computations on the fly in a
post-processor enables extremely accurate field evaluations in the
entire computational domain, also close to surfaces where integral
equation techniques usually encounter difficulties.

Several disparate computational techniques are used. The paper is,
consequently, divided into shorter sections that dwell on specific
issues. Two interleaved overview sections help the reader navigate the
text. The outline is as follows: Sections~\ref{sec:nota}
and~\ref{sec:split} explain our notation and list basic equations.
Section~\ref{sec:Fou} introduces azimuthal Fourier transforms and
presents the integral equation that we actually solve.
Section~\ref{sec:discI} reviews challenges and overall strategies
associated with discretization. Section~\ref{sec:SK} is about a kernel
splitting for integration in the azimuthal direction.
Section~\ref{sec:ZD} provides a link between transformed kernels and
special functions whose evaluation is discussed in
Section~\ref{sec:Q}. Section~\ref{sec:real} presents ideas behind the
kernel-split product integration scheme used in the polar direction in
Sections~\ref{sec:extract} and~\ref{sec:Kastperp}. The entire
discretization scheme is then summarized in Section~\ref{sec:discII}
and illustrated by numerical examples in Section~\ref{sec:numex}.

\begin{figure}
\centering 
\includegraphics[height=47mm]{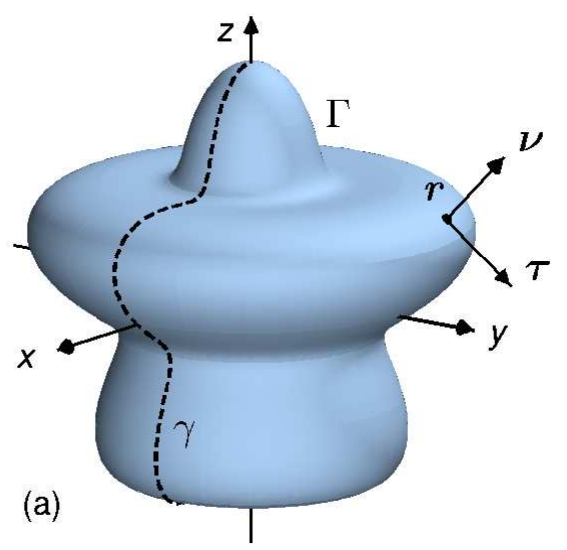}
\includegraphics[height=47mm]{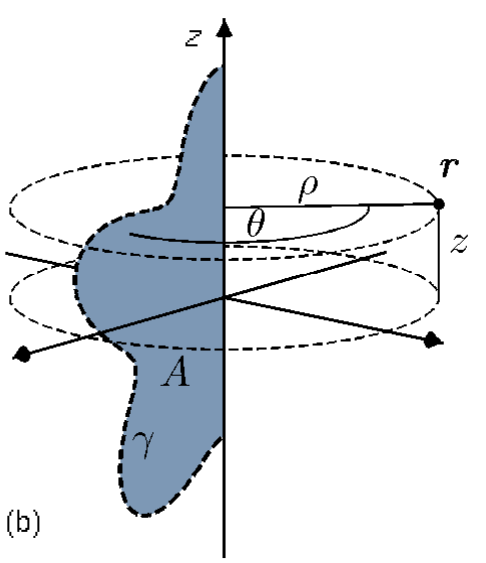}
\includegraphics[height=47mm]{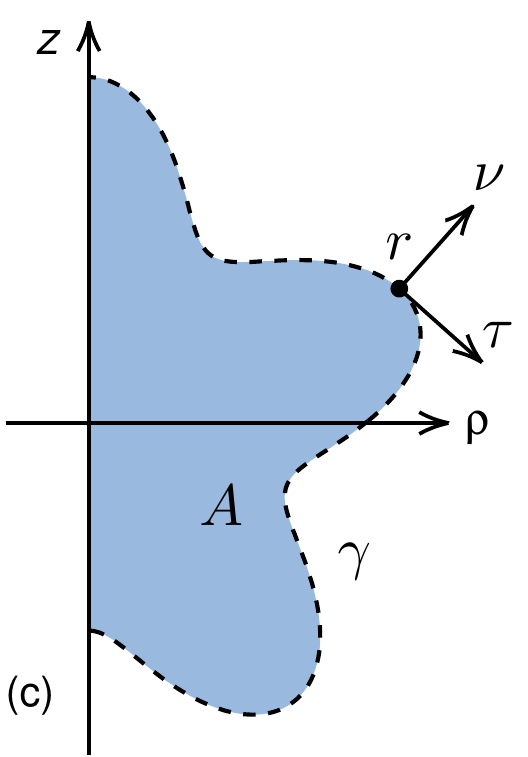}
\caption{\sf An axially symmetric surface $\Gamma$ generated by a curve 
  $\gamma$. (a) A point ${\vec r}$ on $\Gamma$ has outward unit normal
  $\boldsymbol{\nu}$ and tangent vector $\boldsymbol{\tau}$. (b)
  ${\vec r}$ has radial distance $\rho$, azimuthal angle $\theta$, and
  height $z$. The two-dimensional domain $A$ is bounded by $\gamma$
  and the $z$-axis. (c) Coordinate axes and vectors in the half-plane
  $\theta=0$.}
\label{fig:geometry}
\end{figure}

\section{Notation and integral equation} 
\label{sec:nota}

Let $\Gamma$ be an axially symmetric surface enclosing a
three-dimensional domain $V$ (a body of revolution) and let
\begin{displaymath}
{\vec r}=(x,y,z)=(\rho\cos{\theta},\rho\sin{\theta},z)
\end{displaymath}
denote points in $\mathbb{R}^3$. Here $\rho=\sqrt{x^2+y^2}$ is the
radial distance from ${\vec r}$ to the $z$-axis and $\theta$ is the
azimuthal angle. The outward unit normal vector at a point ${\vec r}$
on $\Gamma$ is
$\boldsymbol{\nu}=(\nu_\rho\cos{\theta},\nu_\rho\sin{\theta},\nu_z)$
and
$\boldsymbol{\tau}=(\nu_z\cos{\theta},\nu_z\sin{\theta},-\nu_\rho)$ is a unit tangent vector. See Figure~\ref{fig:geometry}(a)
and~\ref{fig:geometry}(b).

The angle $\theta=0$ defines a half-plane in $\mathbb{R}^3$ whose
intersection with $\Gamma$ corresponds to a generating curve $\gamma$.
We introduce $r=(\rho,z)$ for points in this half-plane and let
$A$ be the two-dimensional closed region bounded by $\gamma$ and the
$z$-axis. The outward unit normal at a point $r$ on $\gamma$ is
$\nu=(\nu_\rho,\nu_z)$ and $\tau=(\nu_z,-\nu_\rho)$ is a
tangent. See Figure~\ref{fig:geometry}(c).

We shall discretize the three layer potential operators $S$,
$K_{\boldsymbol{\nu}}$, and $K_{\boldsymbol{\tau}}$ defined by their
actions on a layer density $\varrho$ on $\Gamma$ as
\begin{align}
S\varrho({\vec r})&=
\int_{\Gamma}
\frac{e^{{\rm i}k|{\vec r}-{\vec r}'|}}
{4\pi|{\vec r}-{\vec r}'|}\varrho({\vec r}')\,{\rm d}\Gamma'\,,
\label{eq:S}\\
K_{\boldsymbol{\nu}}\varrho({\vec r})&=
\int_{\Gamma}\frac{\partial}{\partial\boldsymbol{\nu}}
\frac{e^{{\rm i}k|{\vec r}-{\vec r}'|}}
{4\pi|{\vec r}-{\vec r}'|}\varrho({\vec r}')\,{\rm d}\Gamma'\,,
\label{eq:Kast}\\
K_{\boldsymbol{\tau}}\varrho({\vec r})&=
\int_{\Gamma}\frac{\partial}{\partial\boldsymbol{\tau}}
\frac{e^{{\rm i}k|{\vec r}-{\vec r}'|}}
{4\pi|{\vec r}-{\vec r}'|}\varrho({\vec r}')\,{\rm d}\Gamma'\,.
\label{eq:Kastperp}
\end{align}
Here $k$ is the wavenumber and differentiation with respect to
$\boldsymbol{\nu}$ and to $\boldsymbol{\tau}$ denote normal- and
tangential derivatives. The representation
\begin{equation}
u({\vec r})=S\varrho({\vec r})\,,
\quad {\vec r}\in V\cup\Gamma\,,
\label{eq:field}
\end{equation}
for a solution $u({\vec r})$ to the interior Neumann problem for the
Helmholtz equation with boundary condition $f({\vec r})$ on $\Gamma$
\begin{gather}
  \Delta u({\vec r})+k^2u({\vec r})=0\,, \quad {\vec r}\in V\,,
  \label{eq:Helm1}\\
  \boldsymbol{\nu}\cdot\nabla u({\vec r})=f({\vec r})\,,
  \quad {\vec r}\in\Gamma\,,
  \label{eq:Helm2}
\end{gather}
gives the integral equation
\begin{equation}
\left(I+2K_{\boldsymbol{\nu}}\right)\varrho({\vec r})=2f({\vec r})\,,
\quad {\vec r}\in\Gamma\,.
\label{eq:inteq}
\end{equation}

In order to normalize solutions $u^{\rm h}({\vec r})$
to~(\ref{eq:Helm1}) and~(\ref{eq:Helm2}) with homogeneous boundary
conditions we use
\begin{equation}
\int_V\left|u^{\rm h}({\vec r})\right|^2\,{\rm d}V=
\frac{1}{2k^2}\int_{\Gamma}({\vec r}\cdot\boldsymbol{\nu})
\left(k^2\left|u^{\rm h}({\vec r})\right|^2-\left|\nabla u^{\rm h}({\vec r})
\right|^2\right)\,{\rm d}\Gamma\,,
\label{eq:Barnett}
\end{equation}
which is a special case of a formula due to
Barnett~\cite[eq.~(12)]{Barn06a}. See also Barnett and
Hassell~\cite[eq.~(47)]{Barn14}.

The operators $K_{\boldsymbol{\nu}}$ and $K_{\boldsymbol{\tau}}$ can
be discretized in similar ways. In what follows we concentrate on the
operators $S$ and $K_{\boldsymbol{\nu}}$. We comment on
$K_{\boldsymbol{\tau}}$, needed in~(\ref{eq:Barnett}), only in
situations where its discretization differs from that of
$K_{\boldsymbol{\nu}}$.

\section{Splittings of kernels}
\label{sec:split}

The kernels $S({\vec r},{\vec r}')$ and $K_{\boldsymbol{\nu}}({\vec
  r},{\vec r}')$ of the operators $S$ and $K_{\boldsymbol{\nu}}$ are
weakly singular at ${\vec r}'={\vec r}$. These singularities cause
problems when azimuthal Fourier coefficients of $S({\vec r},{\vec
  r}')$ and $K_{\boldsymbol{\nu}}({\vec r},{\vec r}')$, also called
modal Green's functions or transformed kernels, are to be evaluated
numerically. In this paper we split $S({\vec r},{\vec r}')$ and
$K_{\boldsymbol{\nu}}({\vec r},{\vec r}')$ as
\begin{align}
S({\vec r},{\vec r}')&=
     Z({\vec r},{\vec r}')\left(H_1({\vec r},{\vec r}')
+{\rm i}H_2({\vec r},{\vec r}')\right)\,,
\label{eq:splitting1}\\
K_{\boldsymbol{\nu}}({\vec r},{\vec r}')&=
     D_{\boldsymbol{\nu}}({\vec r},{\vec r}')\left(H_3({\vec r},{\vec r}')
+{\rm i}H_4({\vec r},{\vec r}')\right)\,,
\label{eq:splitting2}
\end{align}
where
\begin{align}
Z({\vec r},{\vec r}')&=\frac{1}{4\pi|{\vec r}-{\vec r}'|}\,,\\
D_{\boldsymbol{\nu}}({\vec r},{\vec r}')&=
 -\frac{\nu_\rho(\rho-\rho'\cos(\theta-\theta'))
    +\nu_z(z-z')}{4\pi|{\vec r}-{\vec r}'|^3}\,,\\
H_1({\vec r},{\vec r}')&=\cos(k|{\vec r}-{\vec r}'|)\,,
\label{eq:EI}\\
H_2({\vec r},{\vec r}')&=\sin(k|{\vec r}-{\vec r}'|)\,,
\label{eq:EII}\\
H_3({\vec r},{\vec r}')&=\cos(k|{\vec r}-{\vec r}'|)
+k|{\vec r}-{\vec r}'|\sin(k|{\vec r}-{\vec r}'|)\,,
\label{eq:FI}\\
H_4({\vec r},{\vec r}')&=\sin(k|{\vec r}-{\vec r}'|)
-k|{\vec r}-{\vec r}'|\cos(k|{\vec r}-{\vec r}'|)\,,
\label{eq:FII}
\end{align}
and 
\begin{equation}
|{\vec r}-{\vec r}'|=\sqrt{\rho^2+\rho'^2
-2\rho\rho'\cos(\theta-\theta')+(z-z')^2}\,.
\end{equation}
Splittings such as~(\ref{eq:splitting1}) and~(\ref{eq:splitting2}) can
facilitate the evaluation of modal Green's functions, as pointed out
in~\cite{Conw10}. See~\cite[Section II]{Vaes12} for a review of other
splitting options for $S({\vec r},{\vec r}')$ and~\cite{Gust10} for
efficient splitting-free modal Green's function evaluation techniques.

\section{Fourier series expansions} 
\label{sec:Fou}

The first step in our discretization scheme for ~(\ref{eq:field})
and~(\ref{eq:inteq}) is an azimuthal Fourier transformation. For this,
several $2\pi$-periodic quantities need to be expanded in Fourier
series. We define the azimuthal Fourier coefficients
\begin{align}
g_n(r)&=\frac{1}{\sqrt{2\pi}}\int_{-\pi}^{\pi}e^{-{\rm i}n\theta}
g({\vec r})\,{\rm d}\theta\,, \quad g=f,\varrho,u\,,\\
G_n(r,r')&=
\frac{1}{\sqrt{2\pi}}\int_{-\pi}^{\pi}e^{-{\rm i}n(\theta-\theta')}
G({\vec r},{\vec r}')\,{\rm d}(\theta-\theta')\,,
\label{eq:GF}
\end{align}
where $G$ can represent any of the functions $D_{\boldsymbol{\nu}}$,
$D_{\boldsymbol{\nu}}H_3$, $D_{\boldsymbol{\nu}}H_4$, $H_1$, $H_3$,
$K_{\boldsymbol{\nu}}$, $D_{\boldsymbol{\tau}}$,
$K_{\boldsymbol{\tau}}$, $S$, $Z$, $ZH_1$, and $ZH_2$. The subscript
$n$ is called the azimuthal index.

Expansion and integration over $\theta'$ gives
for~(\ref{eq:field}) and~(\ref{eq:inteq})
\begin{gather}
u_n(r)=\sqrt{2\pi}\int_{\gamma}S_n(r,r')\varrho_n(r')\rho'
\,{\rm d}\gamma'\,,\quad r\in A\,,
\label{eq:fieldF}\\
\varrho_n(r)+2\sqrt{2\pi}\int_{\gamma}
K_{\nu n}(r,r')\varrho_n(r')\rho'\,{\rm d}\gamma'=2f_n(r)\,,
\quad r\in \gamma\,.
\label{eq:inteqF}
\end{gather}
Solving the full integral equation~(\ref{eq:inteq}) and evaluating the
field $u({\vec r})$ of~(\ref{eq:field}) amounts to solving a series of
modal integral equations~(\ref{eq:inteqF}) for $\varrho_n$, $n=0,\pm
1,\pm 2,\ldots$, and then retrieving $u({\vec r})$ by summation of its
Fourier series. In this paper, since we are chiefly interested in
Neumann Laplace eigenfunctions, we concentrate on
solving~(\ref{eq:inteqF}) and on evaluating $u_n(r)$
of~(\ref{eq:fieldF}) for individual modes $n$. Note that a solution
$\varrho_n(r)$ to~(\ref{eq:inteqF}) corresponds to a modal field $u_n(r)$
whenever $f_n(r)$ is non-zero. The $j$th homogeneous solution
$\varrho_{n,j}(r)$ to~(\ref{eq:inteqF}) corresponds to a modal
eigenfunction $u_{n,j}(r)$ such that
\begin{equation}
u_{n,j}({\vec r})=\frac{1}{\sqrt{2\pi}}u_{n,j}(r)e^{in\theta}
\end{equation}
is a Neumann Laplace eigenfunction of the original Helmholtz
problem~(\ref{eq:Helm1}) and~(\ref{eq:Helm2}). Then $k^2$ is a Neumann
Laplace eigenvalue and $k$ is a Neumann eigenwavenumber $k_{n,j}$.

Expansion for~(\ref{eq:Barnett}) with $u^{\rm h}_n(r)=u_{n,j}(r)$
gives
\begin{multline}
\int_A\left|u_{n,j}(r)\right|^2\rho\,{\rm d}A=\\
\frac{1}{2k_{n,j}^2}\int_{\gamma}(r\cdot\nu)
\left(\left(k_{n,j}^2-\frac{n^2}{\rho^2}\right)\left|u_{n,j}(r)\right|^2
-\left|\tau\cdot\nabla u_{n,j}(r)\right|^2\right)
\rho\,{\rm d}\gamma\,.
\label{eq:BarnettF}
\end{multline}
Modal eigenfunctions, normalized with respect to this energy integral,
are needed in the evaluation of resonances exited by sources, in the
comparison of field strengths for different eigenfunctions, and in
convergence tests~\cite[Chapter 5]{KaKr14}.

The Fourier coefficients of a product of two $2\pi$-periodic functions
$g({\vec r})$ and $h({\vec r})$ with coefficients $g_n(r)$ and
$h_n(r)$ can be obtained by convolution
\begin{equation}
(gh)_n(r)=\frac{1}{\sqrt{2\pi}}\sum_{m=-\infty}^{\infty}g_m(r)h_{n-m}(r)\,.
\label{eq:convolve}
\end{equation}

\section{Discretization -- overview I}
\label{sec:discI}

We seek, for a given $n$, a Nyström discretization
of~(\ref{eq:fieldF}) and of~(\ref{eq:inteqF}). There are two
difficulties here: First, the logarithmically singular kernels
$S_n(r,r')$ and $K_{\nu n}(r,r')$ need to be evaluated at a set of
point-pairs $(r,r')$. Second, suitable quadrature weights need to be
found for integration along $\gamma$.

When $r$ and $r'$ are distant, the kernels $S_n(r,r')$ and
$K_{\nu n}(r,r')$ are evaluated from their definitions~(\ref{eq:GF})
using discrete Fourier transform techniques (FFT) in the azimuthal
direction.

When $r$ and $r'$ are close, we split $S_n(r,r')$ and
$K_{\nu n}(r,r')$ into two parts each: One part which again is
computed directly via FFT and another part which is computed using
convolution of $Z_n(r,r')$ with $H_{1n}(r,r')$ and of
$D_{\nu n}(r,r')$ with $H_{3n}(r,r')$, respectively. These
splittings, originating from~(\ref{eq:splitting1})
and~(\ref{eq:splitting2}), are further discussed in
Section~\ref{sec:SK}. 

The functions $H_{1n}(r,r')$ and $H_{3n}(r,r')$, needed for the
convolution, are computed via FFT. The functions $Z_n(r,r')$ and
$D_{\nu n}(r,r')$, also needed for the convolution, are treated using
semi-analytical techniques described in Section~\ref{sec:ZD}.

The construction of quadrature weights which capture the logarithmic
singularities of $S_n(r,r')$ and $K_{\nu n}(r,r')$, assuming that
$\varrho_n(r)$ is smooth, is described in Sections~\ref{sec:real}
and~\ref{sec:extract}.

For simplicity, all FFT operations are controlled by a single integer
$N$. When doing a Fourier series expansion of a function, $g({\vec
  r})$ say, as to get coefficients
\begin{equation}
g_n(r)=\frac{1}{\sqrt{2\pi}}\int_{-\pi}^{\pi}e^{-{\rm i}n\theta}
g({\vec r})\,{\rm d}\theta\,,
\end{equation}
we use $2N+1$ equispaced points in the azimuthal direction so that
\begin{align}
g_n(r)\approx\frac{\sqrt{2\pi}}{2N+1}\sum_{m=-N}^{N}e^{-{\rm i}n\theta_m}
g(\rho\cos{\theta_m},\rho\sin{\theta_m},z)\,,
\qquad \theta_m=\frac{2\pi m}{2N+1}\,.
\label{eq:trapezo}
\end{align}
We note, but do not exploit, that the need for azimuthal resolution of
$g({\vec r})$ may vary with $z$.

A convolution, such as~(\ref{eq:convolve}), is truncated to
\begin{equation}
(gh)_n(r)=\frac{1}{\sqrt{2\pi}}\sum_{m=\max\{n-N,-N\}}^{\min\{N,N+n\}} 
g_m(r)h_{n-m}(r)\,.
\label{eq:convolve2}
\end{equation}
We note, but do not exploit, that other limits in the sum
of~(\ref{eq:convolve2}) can be more efficient in certain situations.

\section{The transformed kernels $S_n(r,r')$ and $K_{\nu n}(r,r')$} 
\label{sec:SK}

The kernels $S_n(r,r')$ and $K_{\nu n}(r,r')$ of (\ref{eq:GF}),
appearing in~(\ref{eq:fieldF}) and~(\ref{eq:inteqF}), can be split
using~(\ref{eq:splitting1}) and~(\ref{eq:splitting2}) as
\begin{multline}
S_n(r,r')=
 \frac{1}{\sqrt{2\pi}}\int_{-\pi}^{\pi}e^{-{\rm i}n(\theta-\theta')}
Z({\vec r},{\vec r}')H_1({\vec r},{\vec r}')\,{\rm d}(\theta-\theta')\\
+\frac{{\rm i}}{\sqrt{2\pi}}\int_{-\pi}^{\pi}e^{-{\rm i}n(\theta-\theta')}
Z({\vec r},{\vec r}')H_2({\vec r},{\vec r}')\,{\rm d}(\theta-\theta')
\label{eq:Sn}
\end{multline}
and
\begin{multline}
K_{\nu n}(r,r')=
 \frac{1}{\sqrt{2\pi}}\int_{-\pi}^{\pi}e^{-{\rm i}n(\theta-\theta')}
D_{\boldsymbol{\nu}}({\vec r},{\vec r}')H_3({\vec r},{\vec r}')
\,{\rm d}(\theta-\theta')\\
+\frac{{\rm i}}{\sqrt{2\pi}}\int_{-\pi}^{\pi}e^{-{\rm i}n(\theta-\theta')}
D_{\boldsymbol{\nu}}({\vec r},{\vec r}')H_4({\vec r},{\vec r}')
\,{\rm d}(\theta-\theta')\,.
\label{eq:Kastn}
\end{multline}
These splittings are useful when $r$ and $r'$ are close. Then, the
second integrals in~(\ref{eq:Sn}) and~(\ref{eq:Kastn}) have smooth
integrands and are computed via FFT, or using straight-forward
integration if only a single $n$ is of interest. The first integrals
in~(\ref{eq:Sn}) and~(\ref{eq:Kastn}) have non-smooth integrands and
are computed via convolution of $Z_n(r,r')$ with $H_{1n}(r,r')$ and of
$D_{\nu n}(r,r')$ with $H_{3n}(r,r')$. 

We now explain the benefit of this strategy more in detail. For $r$
and $r'$ close, the functions $D_{\boldsymbol{\nu}}H_4$, $H_1$, $H_3$,
$ZH_2$ are smooth while $D_{\boldsymbol{\nu}}$,
$D_{\boldsymbol{\nu}}H_3$, $H_2$, $H_4$, $Z$, $ZH_1$ are non-smooth.
The Fourier coefficients $g_n(r)$ of a smooth function $g({\vec r})$
decay rapidly with $n$ and the individual coefficients $g_n(r)$
converge rapidly with $N$ in the FFT. An individual coefficient
$(gh)_n(r)$ in the convolution~(\ref{eq:convolve2}) of two series of
coefficients $g_m(r)$ and $h_m(r)$ has a rapid asymptotic convergence
with $N$ if at least one series has rapidly decaying coefficients --
compare the discussion of product integration in~\cite[Section
6.1]{Hao14}. If, for $r$ and $r'$ close, no splitting was used and the
second integrals in~(\ref{eq:Sn}) and~(\ref{eq:Kastn}) were to be
computed via convolution along with the first integrals, as
in~\cite{Youn12}, then two slowly decaying series would be convolved
and slower convergence in the azimuthal direction is expected. In
other words, the functions $(ZH_1)_n(r,r')$ and
$(D_{\boldsymbol{\nu}}H_3)_n(r,r')$ converge rapidly with $N$ if
computed via convolution (given that $Z_n(r,r')$ and $D_{\nu n}(r,r')$
are available), but slowly if computed via FFT. The functions
$(ZH_2)_n(r,r')$ and $(D_{\boldsymbol{\nu}}H_4)_n(r,r')$ converge
rapidly with $N$ if computed via FFT, but slowly if computed via
convolution.

A precise definition of what it means that $r$ and $r'$ are close is
given in Section~\ref{sec:quad}.

\section{The functions $Z_n(r,r')$ and $D_{\nu n}(r,r')$} 
\label{sec:ZD}

The functions $Z_n(r,r')$ and $D_{\nu n}(r,r')$ of~(\ref{eq:GF}),
needed for the convolution of the first integrals in~(\ref{eq:Sn})
and~(\ref{eq:Kastn}) when $r$ and $r'$ are close, are evaluated using
semi-analytical techniques and special functions. There are several
ways to proceed. One option is presented in~\cite[Section
III]{Vaes12}. We follow Refs.~\cite{Cohl99,Youn12} and write
\begin{equation}
Z_n(r,r')=\frac{1}{\sqrt{8\pi^3\rho\rho'}}
\mathfrak{Q}_{n-\frac{1}{2}}(\chi)
\label{eq:ZQ}
\end{equation}
and
\begin{equation}
D_{\nu n}(r,r')=\frac{1}{\sqrt{8\pi^3\rho\rho'}}\left[
\left(d_{\nu}(r,r')-\frac{\nu_\rho}{2\rho}\right)
\mathfrak{R}_n(\chi)
-\frac{\nu_\rho}{2\rho}\mathfrak{Q}_{n-\frac{1}{2}}(\chi)\right]\,.
\label{eq:DastRQ}
\end{equation}
Here 
\begin{gather}
\chi=1+\frac{|r-r'|^2}{2\rho\rho'}\,,
\label{eq:chidef}\\
d_{\nu}(r,r')=\frac{\nu\cdot(r-r')}{|r-r'|^2}\,,\\
\mathfrak{R}_n(\chi)=\frac{2n-1}{\chi+1}\left(
\chi\mathfrak{Q}_{n-\frac{1}{2}}(\chi)
-\mathfrak{Q}_{n-\frac{3}{2}}(\chi)\right)\,,\qquad n\ge 0\,,
\label{eq:Rfrak}
\end{gather}
and $\mathfrak{Q}_{n-\frac{1}{2}}(\chi)$ are half-integer degree
Legendre functions of the second kind whose evaluation is discussed in
Section~\ref{sec:Q}. Note that $\chi\ge 1$.

\section{The evaluation of $\mathfrak{Q}_{n-\frac{1}{2}}(\chi)$ and 
$\mathfrak{R}_n(\chi)$}
\label{sec:Q}

The Legendre functions $\mathfrak{Q}_{n-\frac{1}{2}}(\chi)$, which for
$\chi\ge 1$ may be called toroidal harmonics, have logarithmic
singularities at $\chi=1$ but are otherwise analytic. We only need to
consider non-negative integers $n$ since it holds that
$\mathfrak{Q}_{-n-\frac{1}{2}}(\chi)=\mathfrak{Q}_{n-\frac{1}{2}}(\chi)$.
The behavior at infinity is~\cite[eq.~(8.1.3)]{Abra72}
\begin{equation}
\lim_{\chi\to\infty}\mathfrak{Q}_{n-\frac{1}{2}}(\chi)
\propto\chi^{-n-\frac{1}{2}}\,, \qquad n\ge 0\,.
\end{equation}

The functions $\mathfrak{Q}_{n-\frac{1}{2}}(\chi)$ can be evaluated in
several ways. We rely on two methods: forward recursion and backward
recursion. The forward recursion is cheap, but unstable for all
$\chi>1$ and sufficiently high $n$. The backward recursion is stable,
but more expensive. It is particularly expensive for $\chi$ close to
unity.

The forward recursion reads~\cite[eq.~(8.5.3)]{Abra72}
\begin{equation}
\mathfrak{Q}_{n-\frac{1}{2}}(\chi)=
 \frac{4n-4}{2n-1}\chi\mathfrak{Q}_{n-\frac{3}{2}}(\chi)
-\frac{2n-3}{2n-1}\mathfrak{Q}_{n-\frac{5}{2}}(\chi)\,,
\qquad n=2,\ldots,N\,,
\label{eq:forw}
\end{equation}
and is, for $\chi>1$, initiated by~\cite[eqs.~(8.13.3)
and~(8.13.7)]{Abra72}
\begin{align}
\mathfrak{Q}_{-\frac{1}{2}}(\chi)&=
\sqrt{\frac{2}{\chi+1}}K_{\rm cei}\left(\frac{2}{\chi+1}\right)\,,
\label{eq:init1}\\
\mathfrak{Q}_{\frac{1}{2}}(\chi)&=
\chi\sqrt{\frac{2}{\chi+1}}K_{\rm cei}\left(\frac{2}{\chi+1}\right)
-\sqrt{2(\chi+1)}E_{\rm cei}\left(\frac{2}{\chi+1}\right)\,,
\label{eq:init2}
\end{align}
where $K_{\rm cei}(m)$ and $E_{\rm cei}(m)$ are complete elliptic
integrals of the first and second kind, respectively, defined as
\begin{align}
K_{\rm cei}(m)&=
\int_0^{\pi/2}\frac{\,{\rm d}\theta}{\sqrt{1-m\sin^2\theta}}\,,\\
E_{\rm cei}(m)&=\int_0^{\pi/2}\sqrt{1-m\sin^2\theta}\;{\rm d}\theta\,,
\end{align}
Note that the definitions of complete elliptic integrals
in~\cite{Abra72} differ between different sections.

The backward recursion is the forward recursion run backwards. It
starts at step $n=M$ with two randomly chosen values for
$\mathfrak{Q}_{M+\frac{1}{2}}(\chi)$ and
$\mathfrak{Q}_{M+\frac{3}{2}}(\chi)$ and is run down to $n=0$. Then
all function values are normalized so that
$\mathfrak{Q}_{-\frac{1}{2}}(\chi)$ agrees with~(\ref{eq:init1}).
Given that~(\ref{eq:init1}) is accurate to some precision, the values
of $\mathfrak{Q}_{n-\frac{1}{2}}(\chi)$, $n=1,\ldots,N$, have that
same accuracy when $M\gg N$ is sufficiently large. The minimal value
of $M$ which has this property depends on $\chi$ and on $N$,
see~\cite[Section 4.6.1]{Gil07}. Alternatively, for $\chi$ close to
unity, the backward recursion could start at step $n=N$ with
$\mathfrak{Q}_{N-\frac{1}{2}}(\chi)$ computed according
to~\cite{Qexp}. See~\cite[Section~8.15]{Abra72} for a similar
suggestion and~\cite{Gil97} and~\cite[Section~12.3]{Gil07} for even
more options.

In the numerical examples of Section~\ref{sec:numex} we choose
backward recursion with $M=N+80$ for $\chi\ge 1.008$ and forward
recursion for $1<\chi<1.008$. To evaluate~(\ref{eq:init1})
and~(\ref{eq:init2}) we use the {\sc Matlab} function {\tt ellipke},
modified as to give better precision when $\chi$ is close to unity and
$\chi-1$ is known to higher absolute accuracy than $\chi$ itself,
compare~(\ref{eq:chidef}).

The functions $\mathfrak{R}_n(\chi)$ are slightly better behaved than
$\mathfrak{Q}_{n-\frac{1}{2}}(\chi)$ since they are finite at
$\chi=1$. Their values are obtained most easily through their
definition~(\ref{eq:Rfrak}) in terms of
$\mathfrak{Q}_{n-\frac{1}{2}}(\chi)$. One can also use the recursion
\begin{equation}
\mathfrak{R}_n(\chi)=
 \frac{4n-4}{2n-3}\chi\mathfrak{R}_{n-1}(\chi)
-\frac{2n-1}{2n-3}\mathfrak{R}_{n-2}(\chi)\,,
\qquad n=2,\ldots,N\,,
\label{eq:forwR}
\end{equation}
initiated by
\begin{align}
\mathfrak{R}_0(\chi)&=
-\sqrt{\frac{2}{\chi+1}}
E_{\rm cei}\left(\frac{2}{\chi+1}\right)\,,\\
\mathfrak{R}_1(\chi)&=
 \sqrt{\frac{2}{\chi+1}}
\left((\chi-1)K_{\rm cei}\left(\frac{2}{\chi+1}\right)-\chi
E_{\rm cei}\left(\frac{2}{\chi+1}\right)\right)\,.
\end{align}

\section{Product integration for singular integrals}
\label{sec:real}

This section summarizes and extends a high-order accurate panel-based
analytic product integration scheme applicable to integrals whose
kernels contain logarithmic- and Cauchy-type singularities. The scheme
was first presented in~\cite{Hels09} and later adapted to the Nyström
discretization of singular integral operators of planar scattering
theory in~\cite{Hels13,HelsKarl13}. In the present work, the scheme is
used for discretization along $\gamma$ of operators containing the
functions $\mathfrak{Q}_{n-\frac{1}{2}}(\chi)$ and $d_\tau(r,r')$, as
explained in Sections~\ref{sec:extract},~\ref{sec:Kastperp},
and~\ref{sec:discII}. Functions with logarithmic singularities occur
in $Z_n(r,r')$ for $r\in A$ and in $D_{\nu n}(r,r')$ for $r\in\gamma$.
The sum of functions with logarithmic- and Cauchy-type singularities
occurs in $D_{\tau n}(r,r')$.

Consider first the discretization of an integral
\begin{equation}
I_p(r)=\int_{\gamma_p}G(r,r')\varrho(r')\,{\rm d}\gamma'\,,
\label{eq:ex1}
\end{equation}
where $G(r,r')$ is a smooth kernel, $\varrho(r)$ is a smooth layer
density, $\gamma_p$ is a quadrature panel on a curve $\gamma$, and $r$
is a point close to, or on, $\gamma_p$. Let $r(t)=(\rho(t),z(t))$
be a parameterization of $\gamma$. Using $n_{\rm pt}$-point
Gauss--Legendre quadrature with nodes and weights $t_j$ and $w_j$,
$j=1,\ldots,n_{\rm pt}$, on $\gamma_p$ it holds to high accuracy
\begin{equation}
I_p(r)=\sum_j G(r,r_j)\varrho_js_jw_j\,.
\label{eq:GauLeg}
\end{equation}
Here $r_j=r(t_j)$, $\varrho_j=\varrho(r(t_j))$, and $s_j=|{\rm d}r(t_j)/{\rm
  d}t|$. When we discretize~(\ref{eq:fieldF}) and~(\ref{eq:inteqF}) we
shall use a Nyström scheme based on panelwise discretization.

We now proceed to find efficient discretizations for~(\ref{eq:ex1})
when~$G(r,r')$ is not smooth, but can be split and factorized into
smooth parts and parts with known singularities.

\subsection{Logarithmic singularity plus smooth part}
\label{sec:Log}

Consider~(\ref{eq:ex1}) when $G(r,r')$ can be expressed as
\begin{equation}
G(r,r')=\log|r-r'|G^{(1)}(r,r')+G^{(0)}(r,r')\,,
\label{eq:ex2}
\end{equation}
where both $G^{(0)}(r,r')$ and $G^{(1)}(r,r')$ are smooth functions.
Then it holds to high accuracy
\begin{equation}
I_p(r)=\sum_j G^{(0)}(r,r_j)\varrho_js_jw_j
      +\sum_j G^{(1)}(r,r_j)\varrho_js_jw_{{\rm L}j}(r)\,,
\label{eq:fix00}
\end{equation}
where $w_{{\rm L}j}(r)$ are $(n_{\rm pt}-1)$th degree product
integration weights for the logarithmic kernel in~(\ref{eq:ex2}). The
weights $w_{{\rm L}j}(r)$ can be constructed using the analytic method
in~\cite[Section 2.3]{Hels09}.

Adding and subtracting
\begin{displaymath}
\sum_j\log|r-r_j|G^{(1)}(r,r_j)\varrho_js_jw_j
\end{displaymath}
to the right in~(\ref{eq:fix00}), assuming $r\ne r_j$, and
using~(\ref{eq:ex2}) we get the expression
\begin{equation}
I_p(r)=\sum_j G(r,r_j)\varrho_js_jw_j
+\sum_j G^{(1)}(r,r_j)\left[
\frac{w_{{\rm L}j}(r)}{w_j}-\log|r-r_j|\right]\varrho_js_jw_j\,.
\label{eq:fix2b}
\end{equation}
Introducing the logarithmic weight corrections $w_{{\rm L}j}^{\rm
  corr}(r)$ for the terms within square brackets in~(\ref{eq:fix2b})
we arrive at
\begin{equation}
I_p(r)=\sum_j G(r,r_j)\varrho_js_jw_j
+\sum_j G^{(1)}(r,r_j)\varrho_js_jw_jw_{{\rm L}j}^{\rm corr}(r)\,.
\label{eq:fix2}
\end{equation}

The expression~(\ref{eq:fix2}) is, from a strictly mathematical
viewpoint, merely~(\ref{eq:fix00}) rearranged in a more appetizing
form without explicit reference to $G^{(0)}$. An important advantage
of~(\ref{eq:fix2}) over~(\ref{eq:fix00}) is, however, related to
computations and appears whenever $r$ coincides with a discretization
point $r_i$ on $\gamma$. Then the expressions for the weight
corrections simplify greatly. In fact, $w_{{\rm L}j}^{\rm corr}(r_i)$,
$i\ne j$, only depends on the relative length (in parameter) of the
quadrature panels upon which $r_i$ and $r_j$ are situated and on nodes
and weights on a canonical panel. See Appendix~A.

For $r=r_j$ in~(\ref{eq:fix2}), neither $G(r_j,r_j)$ nor $w_{{\rm
    L}j}^{\rm corr}(r_j)$ are defined. We revert to~(\ref{eq:fix00})
and set $w_jw_{{\rm L}j}^{\rm corr}(r_j)=w_{{\rm L}j}(r_j)$ and
$G(r_j,r_j)=G^{(0)}(r_j,r_j)$. Apart from that, no explicit knowledge
of $G^{(0)}(r,r')$ is needed in order to implement~(\ref{eq:fix2}). It
suffices to know $G(r,r')$ and $G^{(1)}(r,r')$ numerically at a set of
points.

\medskip\noindent 
{\bf Remark}: We note, but do not exploit, that if
$G^{(1)}(r,r')=\mathcal{O}(|r-r'|)$ as $r\to r'$ one can factor out
$|r-r'|$ from $G^{(1)}(r,r')$ and construct product integration for
the kernel $|r-r'|\log|r-r'|$ rather than for $\log|r-r'|$.

\subsection{Logarithmic- and Cauchy-type singularities plus smooth part}
\label{sec:LogCau}

Now consider~(\ref{eq:ex1}) when $G(r,r')$ can be expressed as
\begin{equation}
G(r,r')=\log|r-r'|G^{(1)}(r,r')
+\frac{\mu\cdot(r-r')}{|r-r'|^2}G^{(2)}(r,r')+G^{(0)}(r,r')\,,
\label{eq:ex3}
\end{equation}
where $G^{(0)}(r,r')$, $G^{(1)}(r,r')$, and $G^{(2)}(r,r')$ are smooth
and $\mu$ is a unit vector. If $\mu=\nu$ and $r\in\gamma$, then the
second kernel on the right in~(\ref{eq:ex3}) is smooth and we are back
to~(\ref{eq:ex2}). Otherwise we proceed as in Section~\ref{sec:Log}
and observe that it holds to high accuracy
\begin{multline}
I_p(r)=\sum_j G^{(0)}(r,r_j)\varrho_js_jw_j
+\sum_j G^{(1)}(r,r_j)\varrho_js_jw_{{\rm L}j}(r)\\
+\sum_j G^{(2)}(r,r_j)\varrho_jw_{{\rm C}j}(r)\,,
\label{eq:fix000}
\end{multline}
where $w_{{\rm L}j}(r)$ are as in~(\ref{eq:fix00}) and $w_{{\rm
    C}j}(r)$ are $(n_{\rm pt}-1)$th degree product integration weights
for the Cauchy-type kernel in~(\ref{eq:ex3}). The weights $w_{{\rm
    C}j}(r)$ can be constructed using the analytic method
in~\cite[Section 2.1]{Hels09}.

Adding and subtracting
\begin{displaymath}
\sum_j\log|r-r_j|G^{(1)}(r,r_j)\varrho_js_jw_j\quad{\rm and}\quad
\sum_j\frac{\mu\cdot(r-r_j)}{|r-r_j|^2}|G^{(2)}(r,r_j)\varrho_js_jw_j
\end{displaymath}
to the right in~(\ref{eq:fix000}), assuming $r\ne r_j$, and
using~(\ref{eq:ex3}) we arrive at an expression of the form
\begin{multline}
I_p(r)=\sum_j G(r,r_j)\varrho_js_jw_j
+\sum_j G^{(1)}(r,r_j)\varrho_js_jw_jw_{{\rm L}j}^{\rm corr}(r)\\
+\sum_j G^{(2)}(r,r_j)\varrho_jw_{{\rm C}j}^{\rm cmp}(r)\,,
\label{eq:fix3}
\end{multline}
where $w_{{\rm C}j}^{\rm cmp}(r)$ are Cauchy-type singular
compensation weights.

The weights $w_{{\rm C}j}^{\rm cmp}(r_i)$ can, similarly to $w_{{\rm
    L}j}^{\rm corr}(r_i)$, be constructed in a particularly economical
way. See Appendix~B. If $r=r_j$ we revert to~(\ref{eq:fix000}) and set
$G(r_j,r_j)=G^{(0)}(r_j,r_j)$. Apart from that, no explicit knowledge
of $G^{(0)}(r,r')$ is needed in order to implement~(\ref{eq:fix3}). It
suffices to know $G(r,r')$, $G^{(1)}(r,r')$, and $G^{(2)}(r,r')$
numerically at a set of points.

\section{Extracting the singularity of
  $\mathfrak{Q}_{n-\frac{1}{2}}(\chi)$}
\label{sec:extract}

The function $\mathfrak{Q}_{n-\frac{1}{2}}(\chi)$ can be split into a
logarithmically singular part and a remainder whenever
$\chi\in(1,3)$~\cite{Qexp}. The splitting reads
\begin{equation}
\mathfrak{Q}_{n-\frac{1}{2}}(\chi)=
-\frac{1}{2}\log\left(\chi-1\right) 
{}_2F_1\left(-n+\frac{1}{2},n+\frac{1}{2};1;\frac{1-\chi}{2}\right)
+R\left(\chi,n\right)\,.
\label{eq:expand}
\end{equation}
Here $R(\chi,n)$ is smooth, $R(1,n)=\log(2)/2+\psi(1)-\psi(n+1/2)$
where $\psi(x)$ is the digamma function, and ${}_2F_1(a,b;c;x)$ is the
hypergeometric function~\cite[eq.~(15.1.1)]{Abra72}
\begin{equation}
{}_2F_1(a,b;c;x)=\sum_{k=0}^{\infty}\frac{(a)_k(b)_k}{(c)_k}\frac{x^k}{k!}\,,
\qquad |x|<1\,,
\label{eq:hyper}
\end{equation}
where $(\cdot)_k$ is the Pochhammer symbol. See
also~\cite[p.~175]{Magn66} for an alternative expression of
$\mathfrak{Q}_\lambda(\chi)$ compatible with~(\ref{eq:expand}) when
$\lambda$ is a half-integer.

\subsection{Use of the splitting}
\label{sec:use}

The splitting~(\ref{eq:expand}) is useful for the discretization
of~(\ref{eq:fieldF}) and~(\ref{eq:inteqF}). This is so since the
singular nature of the kernels $S_n(r,r')$ and $K_{\nu n}(r,r')$ is
contained in $\mathfrak{Q}_{n-\frac{1}{2}}(\chi)$, see~(\ref{eq:Sn}),
(\ref{eq:Kastn}), (\ref{eq:ZQ}), and (\ref{eq:DastRQ}), and
since~(\ref{eq:expand}) expresses $\mathfrak{Q}_{n-\frac{1}{2}}(\chi)$
as a sum of a smooth function and a product of a smooth function and a
logarithmic kernel, provided $\chi<3$. In Section~\ref{sec:Log} we
reviewed kernel-split product integration techniques for the accurate
panel-wise discretization of operators with kernels of this type.
Comparing~(\ref{eq:expand}) to~(\ref{eq:ex2}) one can see that the
expressions are of the same form with ${}_2F_1$ corresponding to minus
$G^{(1)}$.

Within our overall discretization scheme, the kernels
of~(\ref{eq:fieldF}) and~(\ref{eq:inteqF}) are expressed in terms of
$\mathfrak{Q}_{n-\frac{1}{2}}(\chi)$ only when $r$ and $r'$ are close.
We highlight the dependence on $|r-r'|$
in~$\mathfrak{Q}_{n-\frac{1}{2}}(\chi)$ by
rewriting~(\ref{eq:expand}), using~(\ref{eq:chidef}), as
\begin{equation}
\mathfrak{Q}_{n-\frac{1}{2}}(\chi)=
-\log\left|r-r'\right|{}_2F_1\left(
-n+\frac{1}{2},n+\frac{1}{2};1;-\frac{|r-r'|^2}{4\rho\rho'}\right)
+R\left(r,r',n\right)\,.
\label{eq:expand2}
\end{equation}
Here $R(r,r',n)$ is a new remainder which for $r'=r$ assumes the value
\begin{equation}
R(r,r,n)=\log(2\rho)+\psi(1)-\psi\left(n+\frac{1}{2}\right)\,.
\end{equation}

From a numerical viewpoint, there are some problems
with~(\ref{eq:expand2}). The function ${}_2F_1$ may be costly to
compute for a large number of point pairs $(r,r')$ and indices $n$.
Furthermore, even if $r'$ close to $r$ often means that $\chi$ is
close to unity and that $R$ is smooth, this does not have to be the
case when $r$ is close to the endpoints of $\gamma$. For example, the
fourth argument of ${}_2F_1$, and ${}_2F_1$ itself, may then vary
rapidly with $\rho'$ along an individual quadrature panel. The
requirement $\chi<3$ may even be violated so that~(\ref{eq:expand}) is
no longer valid. Similar problems occur for large $n$ in combination
with too wide panels.

To alleviate some of the problems mentioned we truncate the sum
in~(\ref{eq:hyper}) after four terms and introduce
\begin{equation}
{}_2\tilde{F}_1(a,b;x)=\sum_{k=0}^3\frac{(a)_k(b)_k}{(k!)^2}x^k\,,
\label{eq:hypert}
\end{equation}
and expand $\rho'$ around $\rho$ and truncate that
expansion, too, after four terms. The result is
\begin{equation}
\mathfrak{Q}_{n-\frac{1}{2}}(\chi)=
-\log\left|r-r'\right|{}_2\tilde{F}_1
\left(-n+\frac{1}{2},n+\frac{1}{2};-T(r,r')\right)
+\tilde{R}\left(r,r',n\right)\,,
\label{eq:expand3}
\end{equation}
where
\begin{equation}
T(r,r')=\frac{|r-r'|^2}{4\rho^2}
\sum_{k=0}^3\left(\frac{\rho-\rho'}\rho\right)^k\,.
\end{equation}

The new splitting~(\ref{eq:expand3}) is cheaper to implement
than~(\ref{eq:expand2}) and is slightly better balanced. Note that the
new remainder $\tilde{R}(r,r',n)$ is only $\mathcal{C}^5$-smooth away
from the $z$-axis, but that this appears to be sufficient for our
numerical purposes. See Section~\ref{sec:disc1} for additional
techniques used to ensure that the splitting~(\ref{eq:expand3}) can
produce an accurate discretization of integral operators containing
$\mathfrak{Q}_{n-\frac{1}{2}}(\chi)$ within our product integration
framework.

\section{The kernel $K_{\tau n}(r,r')$}
\label{sec:Kastperp}

The treatment of $K_{\tau n}(r,r')$ closely follows that of
$K_{\nu n}(r,r')$, but with $\nu$ replaced by $\tau$. For example,
equation~(\ref{eq:DastRQ}) becomes
\begin{multline}
D_{\tau n}(r,r')=\frac{1}{\sqrt{8\pi^3\rho\rho'}}\left[
d_{\tau}(r,r')\mathfrak{R}_n(\chi)-\frac{\nu_z}{2\rho}
\left(\mathfrak{R}_n(\chi)+\mathfrak{Q}_{n-\frac{1}{2}}(\chi)\right)\right]\,,
\label{eq:Dastperp}
\end{multline}
where $d_{\tau}(r,r')$ is the Cauchy-type singular kernel
\begin{equation}
d_{\tau}(r,r')=\frac{\tau\cdot(r-r')}{|r-r'|^2}\,.
\end{equation}

When $\chi$ is close to unity we can, for the purpose of
discretization, combine~(\ref{eq:Rfrak}) and~(\ref{eq:expand3}) and
write $d_{\tau}(r,r')\mathfrak{R}_n(\chi)$ in~(\ref{eq:Dastperp}) as
\begin{equation}
d_{\tau}(r,r')\mathfrak{R}_n(\chi)=
\log|r-r'|G^{(1)}(r,r')-d_{\tau}(r,r')+G^{(0)}(r,r')\,,
\label{eq:singkern}
\end{equation}
where
\begin{align}
G^{(0)}(r,r')&=d_{\tau}(r,r')\left(\frac{2n-1}{\chi+1}\left(
\chi\tilde{R}(r,r',n)-\tilde{R}(r,r',n-1)\right)+1\right)\,,\\
\begin{split}
G^{(1)}(r,r')&=
-d_{\tau}(r,r')\frac{2n-1}{\chi+1}\left(\chi
{}_2\tilde{F}_1\left(-n+\frac{1}{2},n+\frac{1}{2};-T(r,r')\right)
\right.\\
&\qquad\qquad\qquad\qquad\quad\left.
-{}_2\tilde{F}_1\left(-n+\frac{3}{2},n-\frac{1}{2};-T(r,r')\right)
\right)\,.
\end{split}
\end{align}
The expression~(\ref{eq:singkern}) is of the type~(\ref{eq:ex3}) with
$\mu=\tau$ and $G^{(2)}(r,r')=-1$. The limits of $G^{(0)}(r,r')$ and
$G^{(1)}(r,r')$ are zero as $r'\to r$.

\section{Discretization -- overview II}
\label{sec:discII}

Our discretization scheme for~(\ref{eq:field}) and~(\ref{eq:inteq})
contains a large number of steps and computational techniques. Now
that most of these have been reviewed, and for ease of reading, we
again summarize the main features of the scheme. We also provide
important implementational details.

\subsection{Quadrature techniques used}
\label{sec:quad}

Several quadrature techniques are involved. In the azimuthal direction
we either use the composite trapezoidal rule or semi-analytical
methods combined with FFT and convolution. In the polar direction,
where the Nyström scheme is applied, we either use composite
Gauss--Legendre quadrature or kernel-split product integration and 16
discretization points per panel. With $n_{\rm pan}$ panels on $\gamma$
and notation as in Section~\ref{sec:real}, eq.~(\ref{eq:inteqF})
assumes the general form
\begin{equation}
\varrho_n(r_i)+2\sqrt{2\pi}\sum_{j=1}^m
K_{\nu n}(r_i,r_j)\varrho_n(r_j)\rho_js_jw_{ij}=2f_n(r_i)\,,
\quad i=1,\ldots,m\,,
\label{eq:inteqFdisc}
\end{equation}
where $m=16n_{\rm pan}$, $\rho_j=\rho(t_j)$, $t_j$ and $w_{ij}$ are
nodes and weights on $\gamma$, and the discretization points play the
role of both target points $r_i$ and source points $r_j$. We say that
the $16n_{\rm pan}$ points $r_i$ constitute a {\it global grid} on
$\gamma$.

The mesh of quadrature panels on $\gamma$ is approximately uniform.
Discretization points located on the same panel or on neighboring
panels are said to be {\it close}. Point pairs that are not close are
said to be {\it distant}. It is the interaction between close point
pairs that may require semi-analytical methods, convolution, and
product integration. Panelwise discretization for distant point pairs
is easy: all kernels are considered smooth and we rely exclusively on
the underlying quadrature, that is, the trapezoidal rule and composite
Gauss--Legendre quadrature with $n_{\rm pt}=16$ and with weights
$w_{ij}$ of~(\ref{eq:inteqFdisc}) independent of $i$. The same
philosophy is used in~\cite{Youn12}.

The discretization~(\ref{eq:inteqFdisc}) is a linear system for
$16n_{\rm pan}$ unknown pointwise values of the layer density
$\varrho_n(r)$. The system matrix can, based on panel affiliation, be
partitioned into $n_{\rm pan}\times n_{\rm pan}$ square blocks with
256 entries each. All interaction between close point pairs is
contained in the block tridiagonal part of this partitioned matrix.

\subsection{The discretization of the integral operator in~(\ref{eq:inteqF})}
\label{sec:disc1}

Let $B$ denote the partitioned matrix corresponding to the system
matrix in~(\ref{eq:inteqFdisc}) and let $B^{(3)}$ denote its block
tridiagonal part. All entries of $B$ that lie outside of $B^{(3)}$ are
evaluated using underlying quadrature. The same holds for
contributions to entries of $B^{(3)}$ that stem from the second
integral in~(\ref{eq:Kastn}). Contributions to entries of $B^{(3)}$
that stem from the first integral in~(\ref{eq:Kastn}) are computed via
convolution of $D_{\nu n}(r,r')$ with $H_{3n}(r,r')$. The functions
$H_{3n}(r,r')$ are computed via FFT. The functions $D_{\nu n}(r,r')$
are computed via~(\ref{eq:DastRQ}) and the evaluation techniques of
Section~\ref{sec:Q}. The quadrature weights associated with $D_{\nu
  n}(r,r')$ are found using the product integration of
Section~\ref{sec:real}.

The discretization of~(\ref{eq:inteqF}) for $r$ and $r'$ both close to
the endpoints of $\gamma$ poses an extra challenge related to the
rapid variation of ${}_2\tilde{F}_1$ and $\tilde{R}$, see
Section~\ref{sec:use}. Therefore we temporarily refine the panels
closest to the $z$-axis by binary subdivision $n_{\rm sub}$ times in
the direction towards the $z$-axis. Then we discretize on this refined
mesh and interpolate the result back to target- and source points on
the global grid. This procedure seems to yield fully accurate results
with $n_{\rm sub}=9$ for $n_{\rm pt}=16$ and affects entries in the
top left $2\times 2$ blocks and bottom right $2\times 2$ blocks of
$B^{(3)}$.

High indices $n$, relative to the spacing of points on $\gamma$, also
require extra care in the product integration for
$\mathfrak{Q}_{n-\frac{1}{2}}(\chi)$. This is so since for large $n$,
the splittings of Section~\ref{sec:extract} are appropriate only in a
narrow zone around $\chi=1$. Away from this zone, the functions
corresponding to $G^{(0)}$ and $G^{(1)}$ in~(\ref{eq:ex2}) behave
badly (become large, converge slowly when expressed as infinite sums,
and suffer from cancellation). The truncation technique of
Section~\ref{sec:use} is not powerful enough to counterbalance this
effect on too wide panels. This problem, again, is remedied with
temporary refinement. Each quadrature panel is temporarily divided
into at most four subpanels with $n_{\rm pt}$ auxiliary discretization
points each. The resulting discretization is then interpolated back to
the global grid. This procedure affects all entries of $B^{(3)}$. Note
that the 10th order accurate Kolm--Rokhlin quadrature
of~\cite{Hao14,Youn12} uses at most 24 auxiliary points per panel and
per target point for a similar purpose.

The integer $N$, controlling the resolution in the azimuthal direction
via~(\ref{eq:trapezo}), is taken to be proportional to $n_{\rm pan}$
with a constant of proportionality depending on the shape of $\gamma$.
Note that, for a fixed $N$ and despite all the special techniques
used, the cost of computing all entries of $B^{(3)}$ grows only
linearly with the number of discretization points on $\gamma$.

\subsection{Improved convergence}
\label{sec:improc}

The quadratures used in~(\ref{eq:inteqFdisc}) have different orders of
accuracy. The trapezoidal rule gives exponential convergence in the
azimuthal direction; composite Gauss--Legendre quadrature with $n_{\rm
  pt}=16$ gives 32nd order accuracy for distant interactions in the
polar direction; product integration with $n_{\rm pt}=16$ gives 16th
order convergence for close interactions.

The convergence of a mixed-quadrature Nyström scheme is controlled by
the error in the quadrature with the lowest order. In our scheme, the
product integration is the weakest link and it is important to make
its error constant small -- something which can be achieved by extra
resolution of known functions in singular kernels. For example, one
can first discretize the parts of~(\ref{eq:inteqF}) that correspond to
entries of $B^{(3)}$ using $n_{\rm pt}=32$ with $n_{\rm sub}=11$ and
then interpolate the result back to the 16 points per panel on the
global grid. See ``scheme B'' of~\cite[Section 8.2]{HelsHols14} for
more details. In the numerical examples below we incorporate a simple
version of this convergence enhancement technique. The result,
typically, is a 30 per cent reduction in the number of global grid
points needed to reach a given accuracy in $\varrho_n$. The number of
kernel evaluations required to form $B^{(3)}$ is, however, increased
by a factor of at least four.

\subsection{The discretization of other integral operators}
\label{sec:disc2}

The discretization of the integral operator in~(\ref{eq:fieldF}) is,
more or less, a subproblem of the discretization of the integral
operator in~(\ref{eq:inteqF}). Field evaluations in a post-processor,
for $r\notin\gamma$, are particularly simple as $r$ and $r'$ never
coincide. There is no need to ``interpolate back to points on the
original panel'' and the product integration weights need not be
stored after use.

The discretization of operators containing the kernel $K_{\tau
  n}(r,r')$ may seem more involved than the discretization of
operators containing $K_{\nu n}(r,r')$. The difference being that the
product integration now involves~(\ref{eq:fix3})
for~(\ref{eq:Dastperp}) rather than~(\ref{eq:fix2})
for~(\ref{eq:DastRQ}). In practice, the extra complication is minor.
The factor corresponding to $G^{(2)}(r,r')$ in~(\ref{eq:fix3}) is a
constant, see~(\ref{eq:singkern}), so the situation is the same as in
the two-dimensional examples treated accurately in~\cite{HelsKarl13}.

\section{Numerical examples}
\label{sec:numex}

Our Fourier--Nyström scheme has been implemented in {\sc Matlab}. We
now test this code for correctness, convergence rate, and achievable
accuracy. The numerical examples cover the determination of modal
fields, boundary value maps, Neumann eigenwavenumbers, and normalized
modal eigenfunctions in the entire computational domain.

Asymptotically, for a body of revolution, the total number of Neumann
Laplace eigenfunctions $u_{n,j}(\vec r)$ with eigenwavenumbers
$k_{n,j}$ bounded by a value $k_{\rm lim}$ is proportional to $k_{\rm
  lim}^3$. For a given index $n$, the number of modal eigenfunctions
$u_{n,j}(r)$ is proportional to $k_{\rm lim}^2$. Eigenwavenumbers with
$n\neq 0$ are always degenerate since, for example, $k_{n,j}$ and
$k_{-n,j}$ are the same with Neumann Laplace eigenfunctions being
complex conjugates to each other. For a given $n$, however, and for
most bodies of revolution, there is no degeneracy. Our examples use
wavenumbers with magnitudes of interest in applications such as
mufflers~\cite{Boij03} and ultrasound spectroscopy~\cite{Ogi02}, where
important wavelengths range from half the diameter of the resonant
volume down to a tenth of the diameter. The magnitudes are also
comparable to the ones used in other numerical work on axisymmetric
Helmholtz problems~\cite{Youn12}.

The code is executed on a workstation equipped with an Intel Core i7
CPU at 3.20 GHz and 64 GB of memory. We refrain from giving extensive
timings since the code is not optimized for execution speed and since
the overall complexity is essentially the same as that of the scheme
in~\cite{Youn12}.

\subsection{Modal field from external point source in domain with
  star shaped cross-section}
\label{sec:point}

This section tests convergence of the solution to a modal interior
Neumann Helmholtz problem~(\ref{eq:fieldF}) and~(\ref{eq:inteqF}) with
$n=1$. The body of revolution $V$, shown in
Figure~\ref{fig:geometry}(a), has a generating curve $\gamma$
parameterized as
\begin{equation}
r(t)=(\rho(t),z(t))
=(1+0.25\cos(5t))(\sin(t),\cos(t))\,,\qquad 0\leq t\leq\pi\,.
\label{eq:star}
\end{equation}
The boundary condition on $\gamma$ is given by the normal derivative
of a field $u_{\rm p}(\vec r)$, excited by a point source at ${\vec
  r}_{\rm p}$ outside $V$
\begin{equation}
u_{\rm p}(\vec r)=5S({\vec r},{\vec r}_{\rm p})=
\frac{5e^{{\rm i}k|{\vec r}-{\vec r}_{\rm p}|}}
{4\pi|{\vec r}-{\vec r}_{\rm p}|}\,,
\qquad k=19,\qquad {\vec r}_{\rm p}=(0.5,0,1)\,.
\label{eq:ps}
\end{equation}
We remark, in connection with~(\ref{eq:star}), that an arc length
parameterization is probably more efficient in terms of resolution.

It follows from~(\ref{eq:ps}) and the definitions in
Sections~\ref{sec:nota} and~\ref{sec:Fou} that the excited modal
fields in $A$ and their derivatives on $\gamma$ are
\begin{align}
  u_{{\rm p}n}(r)&=5S_n(r,r_{\rm p})\,,\qquad r\in A\,,
\label{eq:upn}\\
\nu\cdot\nabla u_{{\rm p}n}(r)&=5K_{\nu n}(r,r_{\rm p})\,,
\qquad r\in \gamma\,,\\
\tau\cdot\nabla u_{{\rm p}n}(r)&=
5K_{\tau n}(r,r_{\rm p})\,,\qquad r\in \gamma\,.
\label{eq:upt}
\end{align}
The source strength is chosen to be five so that the modulus of the
field $u_{{\rm p}1}(r)$ in $V$ peaks approximately at unity. The
wavenumber $k=19$ corresponds to about $7.3$ wavelengths across the
generalized diameter of $V$.

Our scheme determines the modal field $u_1(r)$ by first
solving~(\ref{eq:inteqFdisc}) for $\varrho_1(r_i)$ with $f_1(r_i)=5K_{\nu
  1}(r_i,r_{\rm p})$ and then evaluating a discretization
of~(\ref{eq:fieldF}). The values of $K_{\nu 1}(r_i,r_{\rm p})$ are
obtained via~(\ref{eq:GF}) with $G=K_{\boldsymbol{\nu}}$ and the
trapezoidal rule. The mesh is uniformly refined in parameter $t$ (not
in arc length) with $n_{\rm pan}$ panels corresponding to $16n_{\rm
  pan}$ discretization points on $\gamma$. The integer $N$,
controlling the resolution in the azimuthal direction, is chosen as
$N=4n_{\rm pan}$.

\begin{figure}[t]
\centering 
  \includegraphics[height=51mm]{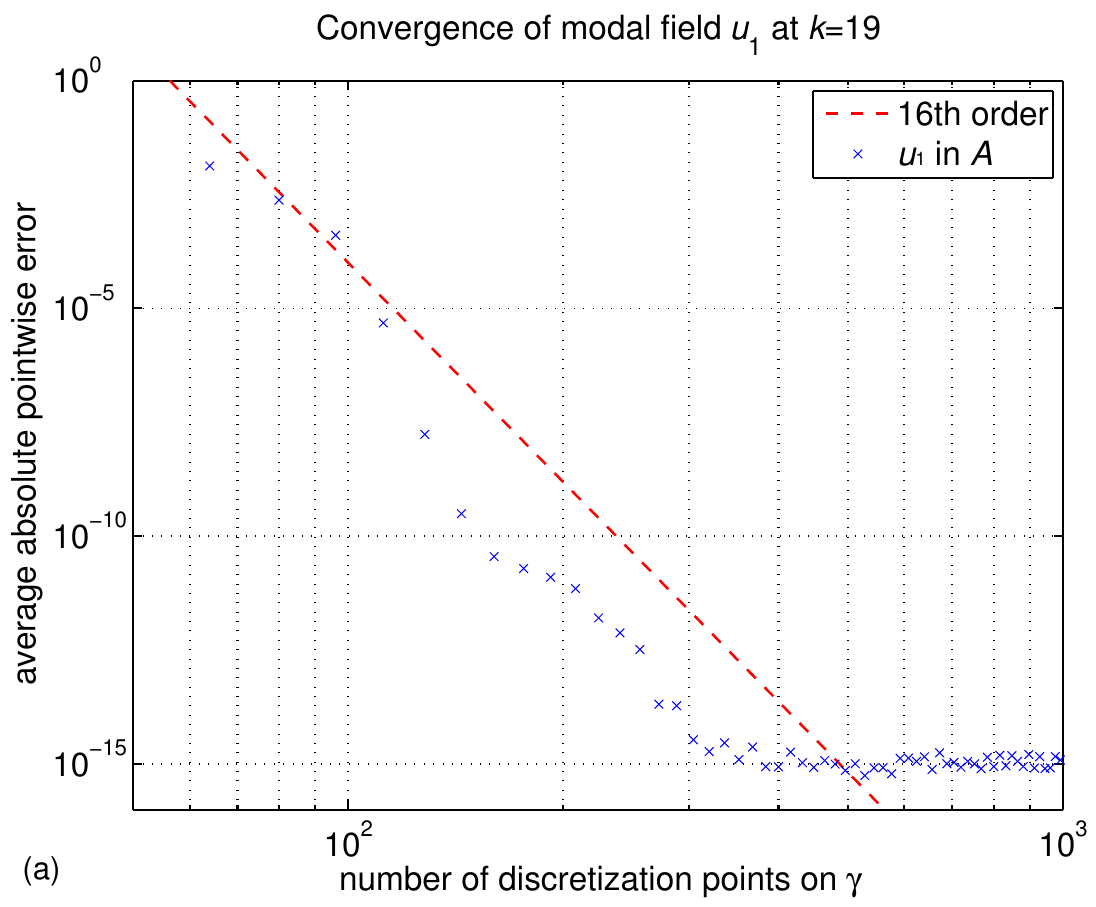}
  \includegraphics[height=51mm]{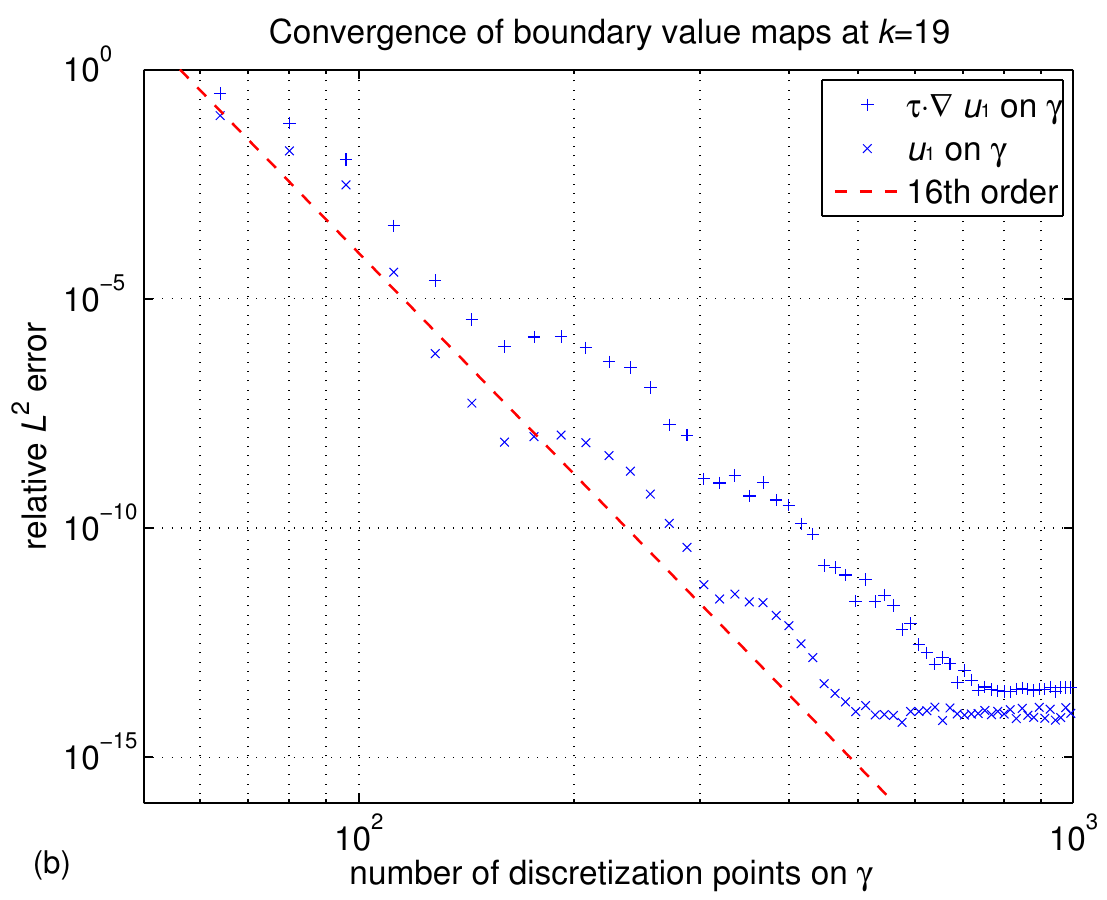}
  \includegraphics[height=50mm]{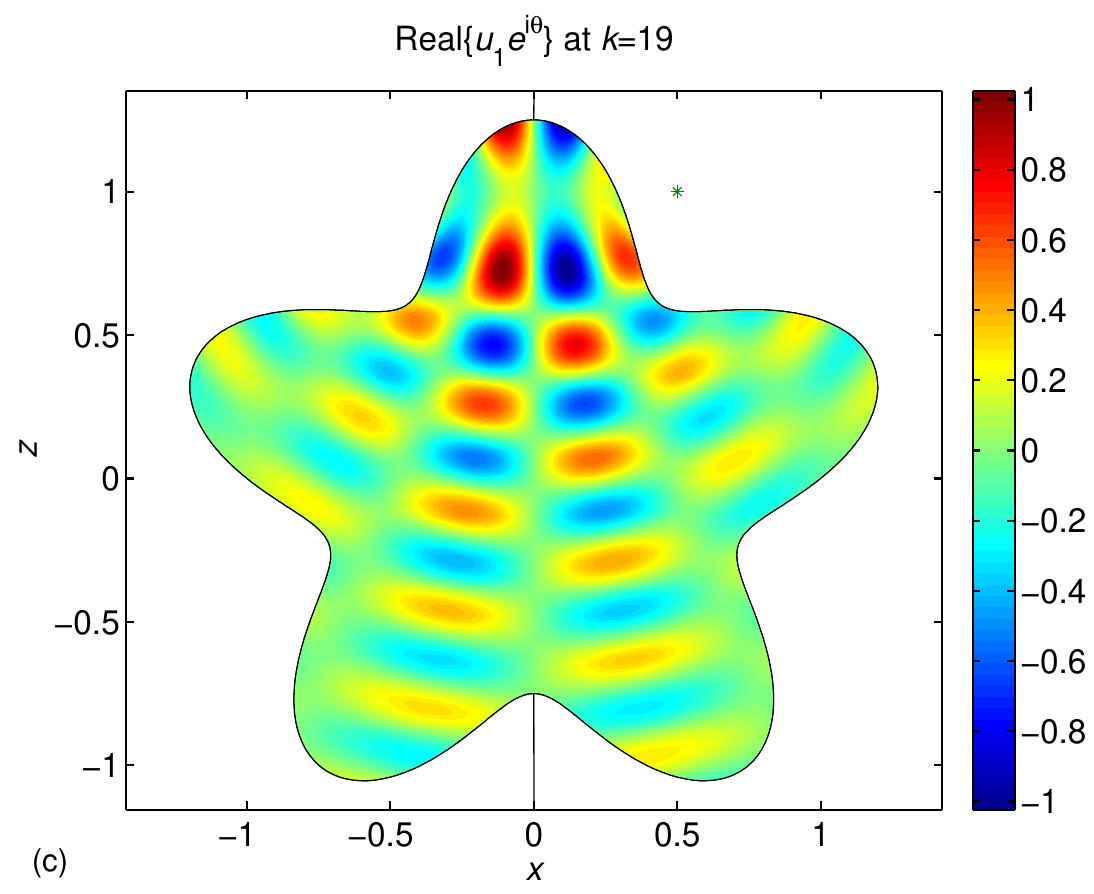}
  \includegraphics[height=50mm]{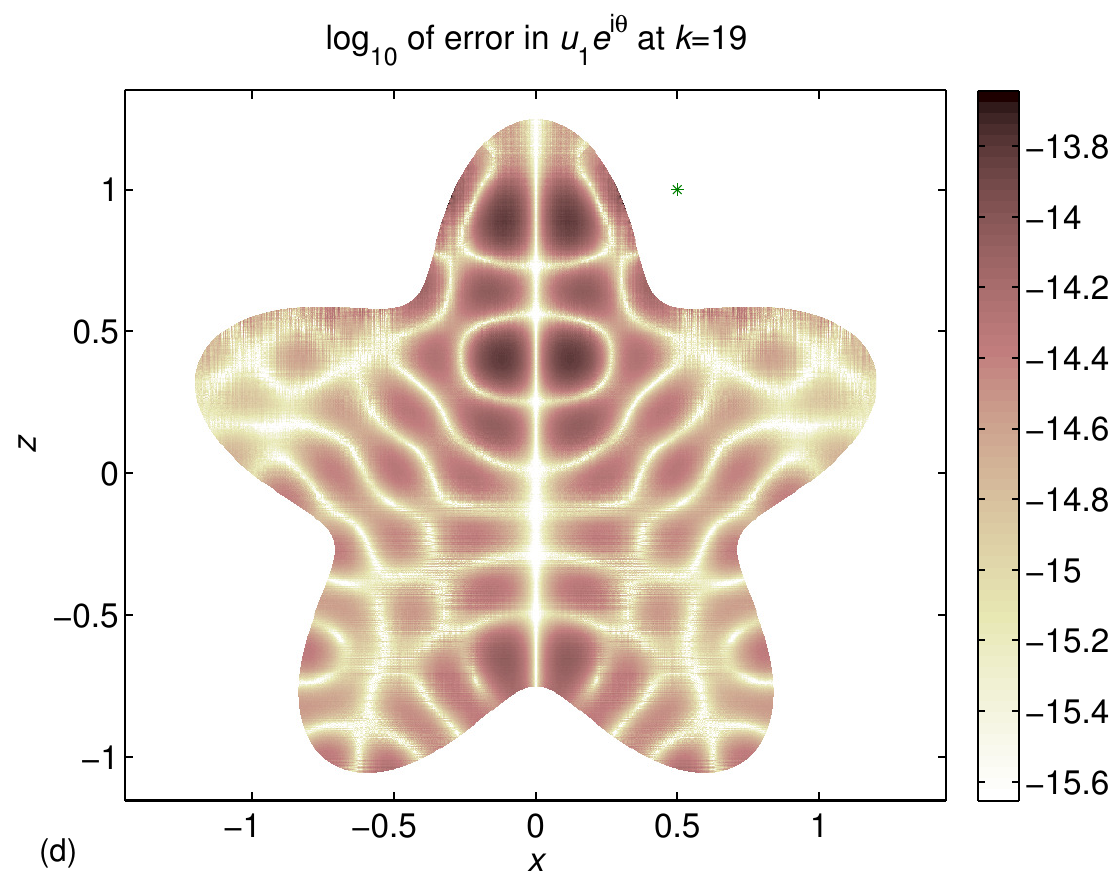}
\caption{\sf Convergence of the modal field $u_1(r)$ excited by a
  point source of strength 5 located at ${\vec r}_{\rm p}=(0.5,0,1)$
  and indicated by a green star. The wavenumber is $k=19$. (a) Average
  absolute pointwise error in $u_1(r)$. (b) Relative $L^2$ errors in
  the modal Neumann-to-Dirichlet map $\nu\cdot\nabla u_1(r)\mapsto
  u_1(r)$ and in the map $\nu\cdot\nabla u_1(r)\mapsto \tau\cdot\nabla
  u_1(r)$. (c) Real part of the field $u_1(r)e^{{\rm i}\theta}$ for
  $\theta=0$ and $\theta=\pi$. (d) $\log_{10}$ of pointwise error in
  $u_1(r)e^{{\rm i}\theta}$ for $\theta=0$ and $\theta=\pi$.}
\label{fig:case1}
\end{figure}

For error estimates we do a comparison with supposedly more accurate
reference values derived directly from~(\ref{eq:ps}).
Figure~\ref{fig:case1}(a) compares our results $u_1(r)$ to reference
values $u_{{\rm p}1}(r)$, obtained from~(\ref{eq:upn})
via~(\ref{eq:GF}) with $G=S$ and the trapezoidal rule. The comparison
is done under mesh refinement and at 50,276 field points in a
cross-section given by the intersection of $V$ and the half-planes
$\theta=0$ and $\theta=\pi$. The field points are placed on a uniform
$300\times 300$ grid in the square $x\in[-1.2,1.2]$ and
$z\in[-1.1,1.3]$. Points outside the cross-section are excluded.
Figure~\ref{fig:case1}(a) shows that the convergence is at least 16th
order, as expected, and that the achievable average accuracy is around
$10\epsilon_{\rm mach}$. The field $u_1(r)$ and the distribution of
absolute pointwise error are depicted in Figures~\ref{fig:case1}(c)
and~\ref{fig:case1}(d). Here the field resolution is increased and
274,800 field points on a uniform $700\times 700$ grid are used. There
are 608 discretization point on $\gamma$. It is worth mentioning that
even though some field points lie very close to $\gamma$, there is no
visible sign of accuracy degradation in the near-surface field
evaluation. We get close to machine precision in the entire
computational domain -- a success which in part can be explained by
the low condition number of the system matrix
in~(\ref{eq:inteqFdisc}). In this example it is only 114.

Having solved~(\ref{eq:inteqFdisc}) for $\varrho_1(r_i)$, we also
evaluate $\tau\cdot\nabla u_1(r)$ at the discretization points $r_i$
on $\gamma$ and compare with reference values. Our values
$\tau\cdot\nabla u_1(r_i)$ are obtained via a discretization of
\begin{equation}
\tau\cdot\nabla u_n(r)=\sqrt{2\pi}\int_{\gamma}
K_{\tau n}(r,r')\varrho_n(r')\rho'
\,{\rm d}\gamma'\,,\quad r\in \gamma\,,
\label{eq:NeuTanF}
\end{equation}
and the techniques of Section~\ref{sec:disc2}. The reference values
$\tau\cdot\nabla u_{{\rm p}1}(r_i)$ are obtained from~(\ref{eq:upt})
via~(\ref{eq:GF}) with $G=K_{\boldsymbol{\tau}}$ and the trapezoidal
rule. Accurate computation of $\tau\cdot\nabla u_n(r_i)$ is important
whenever~(\ref{eq:BarnettF}) is to be used for normalization.
Figure~\ref{fig:case1}(b) shows that the achievable $L^2$ accuracy in
$\tau\cdot\nabla u_1(r)$ is roughly the same as that of $u_1(r)$ on
$\gamma$, albeit somewhat delayed. Had we used numerical
differentiation of $u_1(r)$ for $\tau\cdot\nabla u_1(r)$, rather than
the analytical differentiation implicit in~(\ref{eq:NeuTanF}),
precision would have been lost.

As for timings we quote the following: with $640$ global
discretization points $r_i$ on $\gamma$, and with $2N+1=321$ Fourier
coefficients in the convolutions, it takes 65 seconds in total to
construct the discretizations of $S_1$, $K_{\nu 1}$, and $K_{\tau 1}$,
form and solve~(\ref{eq:inteqFdisc}), and compute $u_1(r_i)$ and
$\tau\cdot\nabla u_1(r_i)$. Of this time, 30 seconds are spent
constructing top left and bottom right blocks of matrices
corresponding to discretizations of $Z_n$, $D_{\nu n}$, and $D_{\tau
  n}$, $n=-N,\ldots,N$, using the procedure for near-endpoint
evaluation described in Sections~\ref{sec:disc1} and~\ref{sec:improc},
and 20 seconds are spent on the remaining tridiagonal blocks of these
$k$-independent matrices.

\subsection{Eigenpair in the unit sphere}
\label{sec:sphere}

This section finds a Neumann eigenwavenumber-eigenfunction pair in the
unit sphere. The generating curve $\gamma$ is parameterized by
\begin{equation}
r(t)=(\rho(t),z(t))=(\sin(t),\cos(t))\,,\qquad 0\leq t\leq\pi\,.
\end{equation}
Eigenpairs in the sphere can be determined from~(\ref{eq:Helm1})
and~(\ref{eq:Helm2}) with $f(\vec r)=0$ using separation of
variables~\cite[Section 9.3]{Arfk05}. Each azimuthal mode $n$ has
infinitely many eigenwavenumbers which can be ordered with respect to
magnitude by two integers $\ell\geq\vert n\vert$ and $m\geq 1$. The
number $k_{n,\ell,m}$ is the $m$th positive solution to
\begin{equation} 
\frac{{\rm d} j_\ell(x)}{{\rm d}x}=0\,,
\label{eq:bessel}
\end{equation} 
where $j_\ell(x)$ is the spherical Bessel function of order $\ell$.

Our approach is to choose a mode $n$, a mesh on $\gamma$, and an
interval $[k_{\rm low},k_{\rm up}]$ where to look for an
eigenwavenumber. For simplicity, we use golden section search with the
condition number of the system matrix in~(\ref{eq:inteqFdisc}) as the
function to be maximized~\cite[Appendix B]{Barn14}. When a maximum is
found, the corresponding $k$ is an approximation of a $k_{n,\ell,m}$.
The integers $\ell$ and $m$ are determined by visual inspection of the
associated eigenfunction. See~\cite{Barn14} for a far more economical
way to determine eigenpairs of the Laplacian using integral equation
techniques.

\begin{figure}[t]
\centering
\includegraphics[height=80mm]{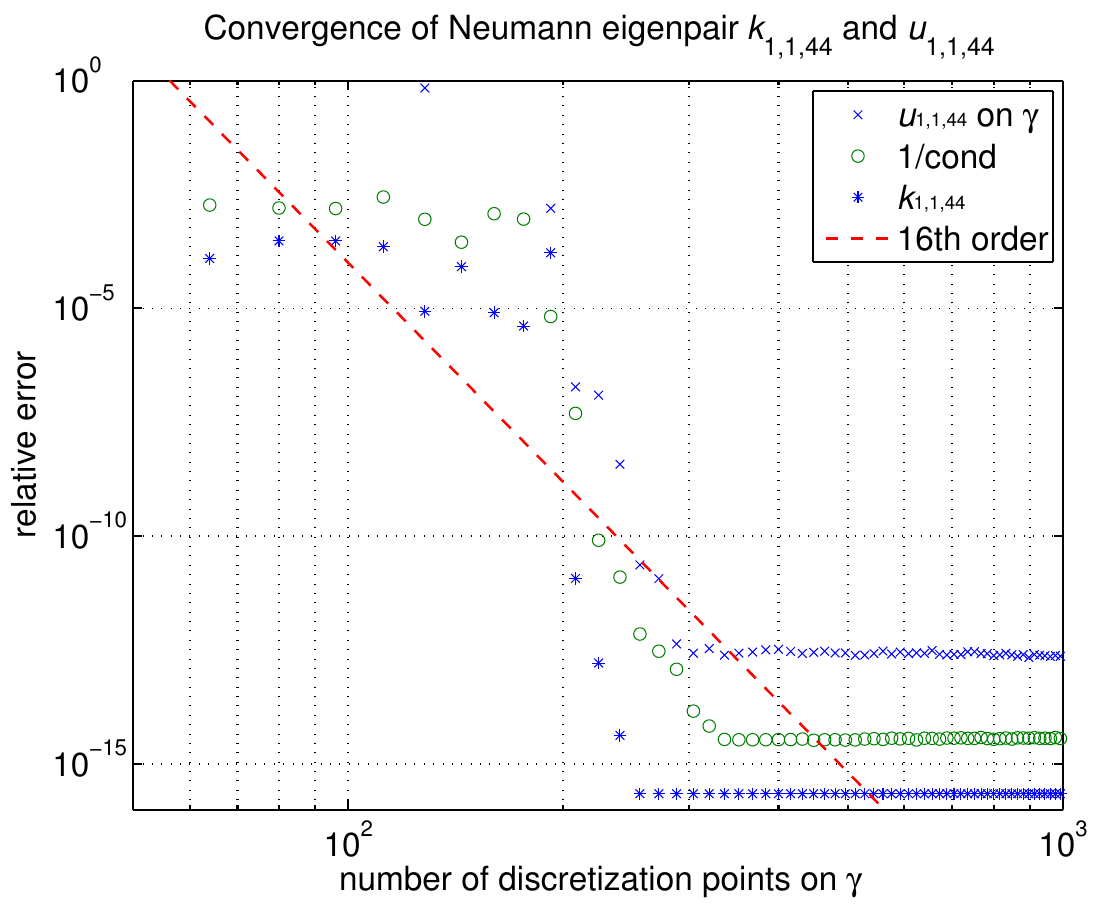}
\caption{\sf Convergence of the Neumann eigenpair $k_{1,1,44}$ and 
  $u_{1,1,44}(r)$ in the unit sphere. The wavenumber
  $k_{1,1,44}=141.3575204174371$ corresponds to approximately 45
  wavelengths across the sphere diameter. x-marks and stars show
  relative errors as a function of mesh resolution. The error in
  $u_{1,1,44}(r)$ is measured in $L^2$ norm at $\gamma$. Circles show
  the reciprocal condition number of the system matrix
  in~(\ref{eq:inteqFdisc}).}
\label{fig:case0}
\end{figure}

The convergence of $k_{1,1,44}$ is shown in Figure~\ref{fig:case0}.
Our computed estimates are compared to the correctly rounded value
$k_{1,1,44}=141.3575204174371$, obtained via~(\ref{eq:bessel}). The
eigenwavenumbers are densely packed so the interval $[k_{\rm
  low},k_{\rm up}]$ has to be narrow. We choose $k_{\rm low}=141.34$
and $k_{\rm up}=141.40$. The mesh is uniformly refined with $n_{\rm
  pan}$ panels corresponding to $16n_{\rm pan}$ discretization points
on $\gamma$. The integer $N$, controlling the resolution in the
azimuthal direction, is chosen as $N=12n_{\rm pan}$.  One can see, in
Figure~\ref{fig:case0}, that a condition number of $10^{12}$ and $256$
points on $\gamma$, corresponding to 3.6 points per wavelength along
$\gamma$, is sufficient to yield $k_{1,1,44}$ with full machine
precision.

Figure~\ref{fig:case0} also shows convergence of the eigenfunction
$u_{1,1,44}(r)$, normalized by
$||u_{n,\ell,m}(r)\sqrt{\rho}||_{L^2(A)}=1$ and by the requirement
that $u_{n,\ell,m}(r)$ is real. A complex constant of unit modulus
which rotates the appropriate eigenvector of the system matrix
in~(\ref{eq:inteqFdisc}), so that it produces a real valued
$u_{n,\ell,m}(r)$ in the directization of~(\ref{eq:fieldF}), is
determined with a least squares fit on $\gamma$. Our computed
estimates for $u_{1,1,44}(r)$ are compared to reference values
obtained from
\begin{equation}
u_{1,1,m}(r(t))=\frac{\sqrt{3}k_{1,1,m}}
{\sqrt{2k_{1,1,m}^2-4}}\sin(t)\,,
\quad r\in\gamma\,.
\end{equation}
The relative error is computed in $L^2$ norm on $\gamma$. The
determination of $u_{1,1,44}(r)$ is a more difficult problem than the
determination of $k_{1,1,44}$. Higher resolution is needed for a given
relative accuracy and the achievable accuracy is lower. 

Even though test problems in the unit sphere often are simple to
solve, the problem in this section may be thought of as a little
harder. The sphere has 204,646 modal eigenfunctions with $k\leq
k_{1,1,44}$ and 2,488 of these are $n=1$ modes. This means that
$k=k_{1,1,44}$, corresponding to approximately 45 wavelengths across
the sphere diameter, is a rather ``high'' wavenumber. See the
interesting discussion in~\cite[Section 1]{Barn14} on what wavenumbers
are needed in applications and on performance characteristics of
different classes of methods used to find them.

\subsection{Eigenpairs in domain with star shaped cross-section}

This section finds two Neumann eigenwavenumbers and their associated
modal eigenfunctions in the body of revolution generated by $\gamma$
of~(\ref{eq:star}). For each azimuthal index $n$ there are infinitely
many Neumann eigenwavenumbers. These can be ordered and numbered with
respect to their magnitude so that $k_{n,j}$ is the $j$th smallest
eigenwavenumber for index $n$. The associated modal eigenfunction
$u_{n,j}(r)$ is normalized, as in Section~\ref{sec:sphere}, by
$||u_{n,j}(r)\sqrt{\rho}||_{L^2(A)}=1$ and by the requirement that
$u_{n,j}(r)$ is real.

\begin{figure}[t]
\centering
\includegraphics[height=51mm]{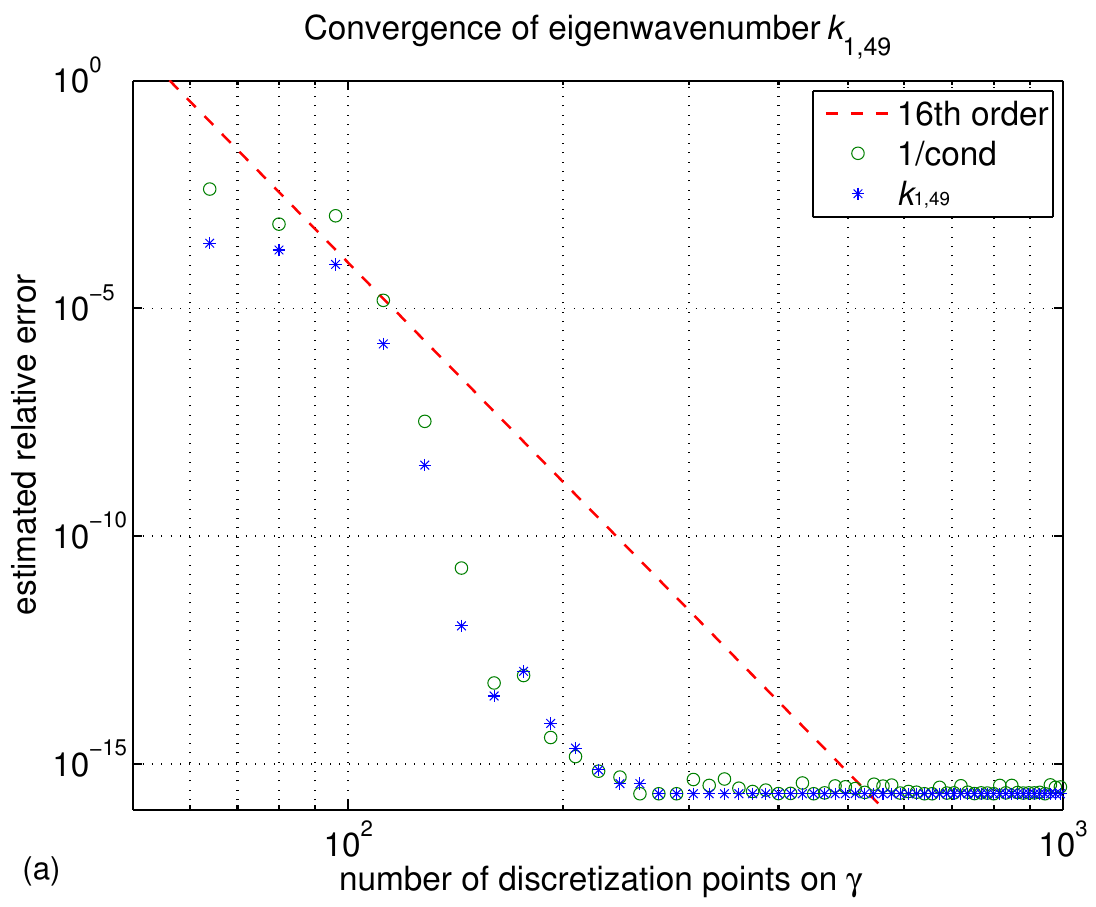}
\includegraphics[height=51mm]{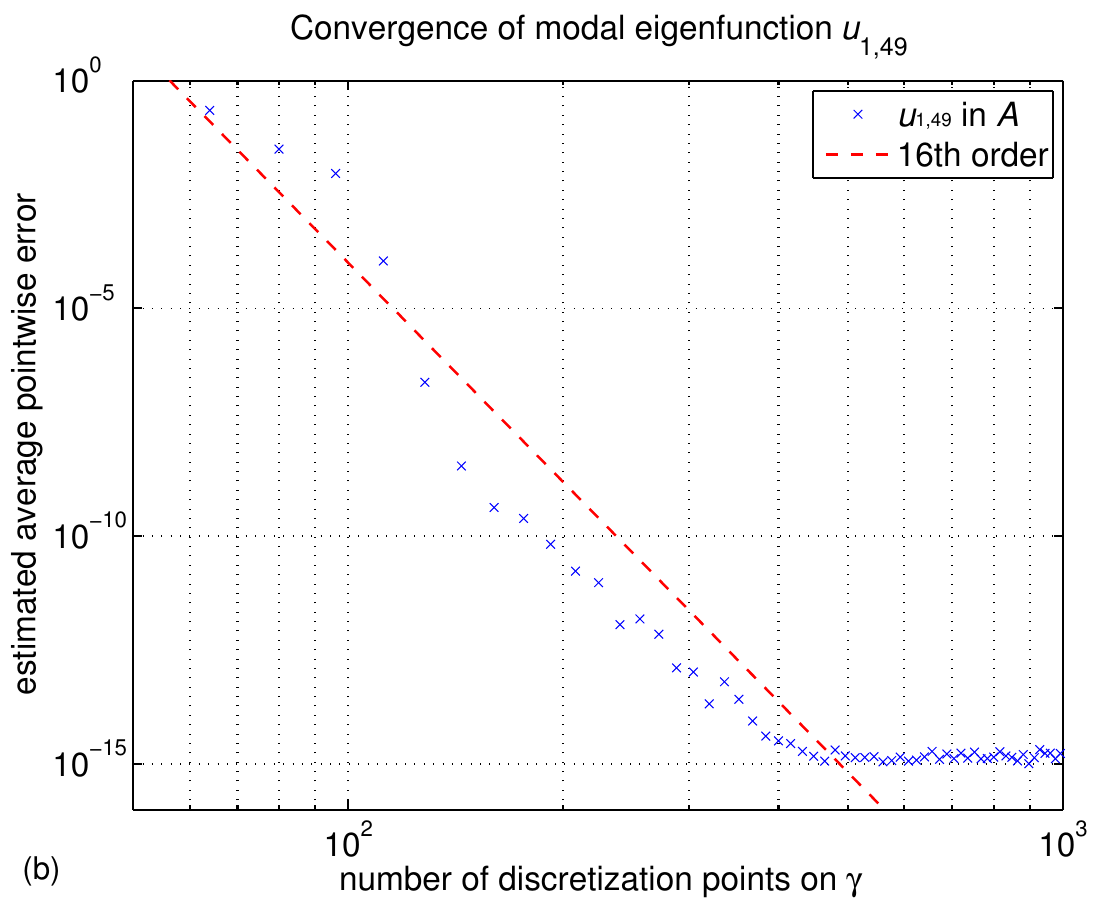}
\includegraphics[height=50mm]{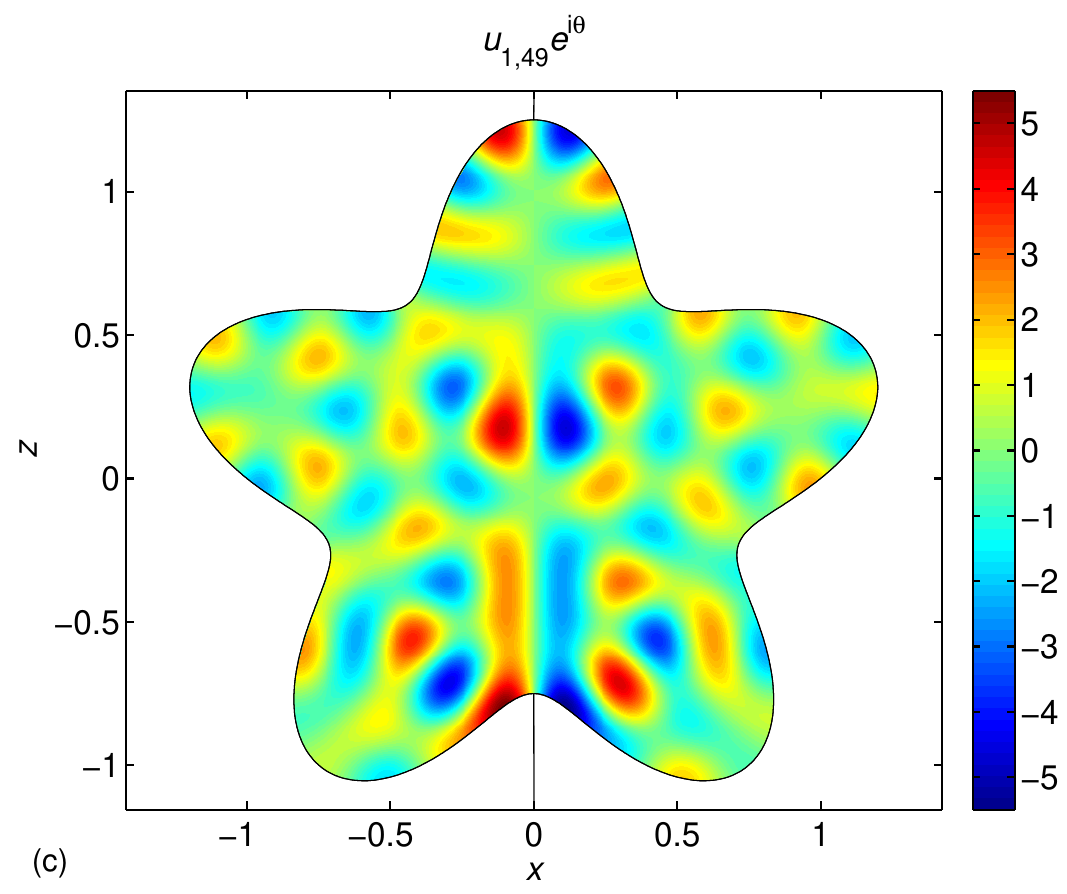}
\includegraphics[height=50mm]{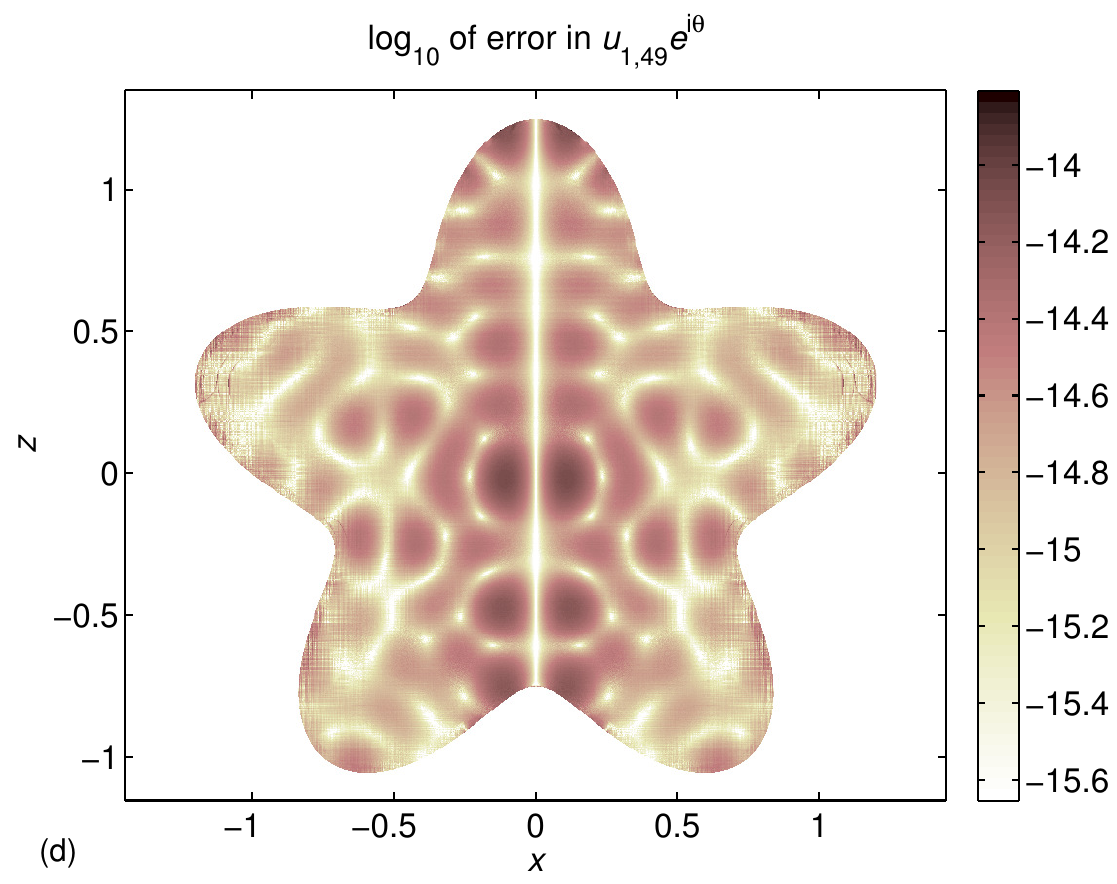}
\caption{\sf Convergence of the Neumann eigenpair $k_{1,49}$ and 
  $u_{1,49}(r)$. (a) Reciprocal condition number and error in
  $k_{1,49}$. Reference value $k_{1,49}=19.22942004015467$ is used.
  (b) Estimated average pointwise error in $u_{1,49}(r)$. (c) The
  field $u_{1,49}(r)e^{{\rm i}\theta}$ for $\theta=0$ and
  $\theta=\pi$. (d) $\log_{10}$ of pointwise error in
  $u_{1,49}(r)e^{{\rm i}\theta}$ for $\theta=0$ and $\theta=\pi$ with
  $560$ points on $\gamma$.}
\label{fig:n1}
\end{figure}

We first search for $k_{1,49}$ in the interval $k_{\rm low}=19.1$ and
$k_{\rm up}=19.3$ using the approach of Section~\ref{sec:sphere}, but
with $N=4n_{\rm pan}$ as in Section~\ref{sec:point}.
Figure~\ref{fig:n1}(a) shows convergence under mesh refinement. The
eigenwavenumber has converged stably to machine precision with 272
discretization points on $\gamma$, corresponding to an average number
of 21.4 points per wavelength along $\gamma$, and we use the converged
value $k_{1,49}=19.229420040154672$ as reference value when estimating
the error.

Having established $k_{1,49}$, we proceed with a convergence study for
the field $u_{1,49}(r)$ in $A$. This investigation is very similar to
the study of the modal field $u_1(r)$ in Section~\ref{sec:point}. The
chief difference being the additional complication of normalizing
$u_{1,49}(r)$ using~(\ref{eq:BarnettF}). Results are displayed in
Figures~\ref{fig:n1}(b), \ref{fig:n1}(c), and~\ref{fig:n1}(d). Here
the {\it pointwise error} refers to an estimated absolute pointwise
error normalized with the largest value of $|u_{1,49}(r)|$, $r\in A$.
The estimated pointwise error at a field point $r$ and with a given
number of discretization points on $\gamma$ is taken as the difference
between the computed value at $r$ and a better resolved value at $r$,
computed with approximately 50 per cent more discretization points on
$\gamma$. Figures~\ref{fig:n1}(a) and~\ref{fig:n1}(b) show asymptotic
16th order convergence. The very high achievable accuracy for the
normalized modal eigenfunction field, shown in Figures~\ref{fig:n1}(b)
and~\ref{fig:n1}(d), is only possible thanks to the accurate
implementation of the boundary value maps $\nu\cdot\nabla
u_n(r)\mapsto u_n(r)$ and $\nu\cdot\nabla u_n(r)\mapsto
\tau\cdot\nabla u_n(r)$, tested separately in
Figure~\ref{fig:case1}(b) and now used in the discretization
of~(\ref{eq:BarnettF}).

As an independent test of correctness in the results for $k_{1,49}$
and $u_{1,49}(r)$ we compared our values and fields with those
obtained from the finite element 2D axisymmetric solver in COMSOL
Multiphysics 4.3b, run on a workstation with 64 GB of memory. With
6,197,297 degrees of freedom in the mesh on $A$, corresponding to 640
degrees of freedom per wavelength, the COMSOL estimates exhibit a
relative deviation from our converged result of about $10^{-10}$ in
$k_{1,49}$ and of $4\cdot 10^{-8}$ in the maximum value of
$|u_{1,49}(r)|$, which for this eigenfunction occurs at $\gamma$. Our
scheme needs roughly 11 points per wavelength for that same accuracy,
see Figure~\ref{fig:n1}(a) and (b).

\begin{figure}[t!]
\centering
\includegraphics[height=50mm]{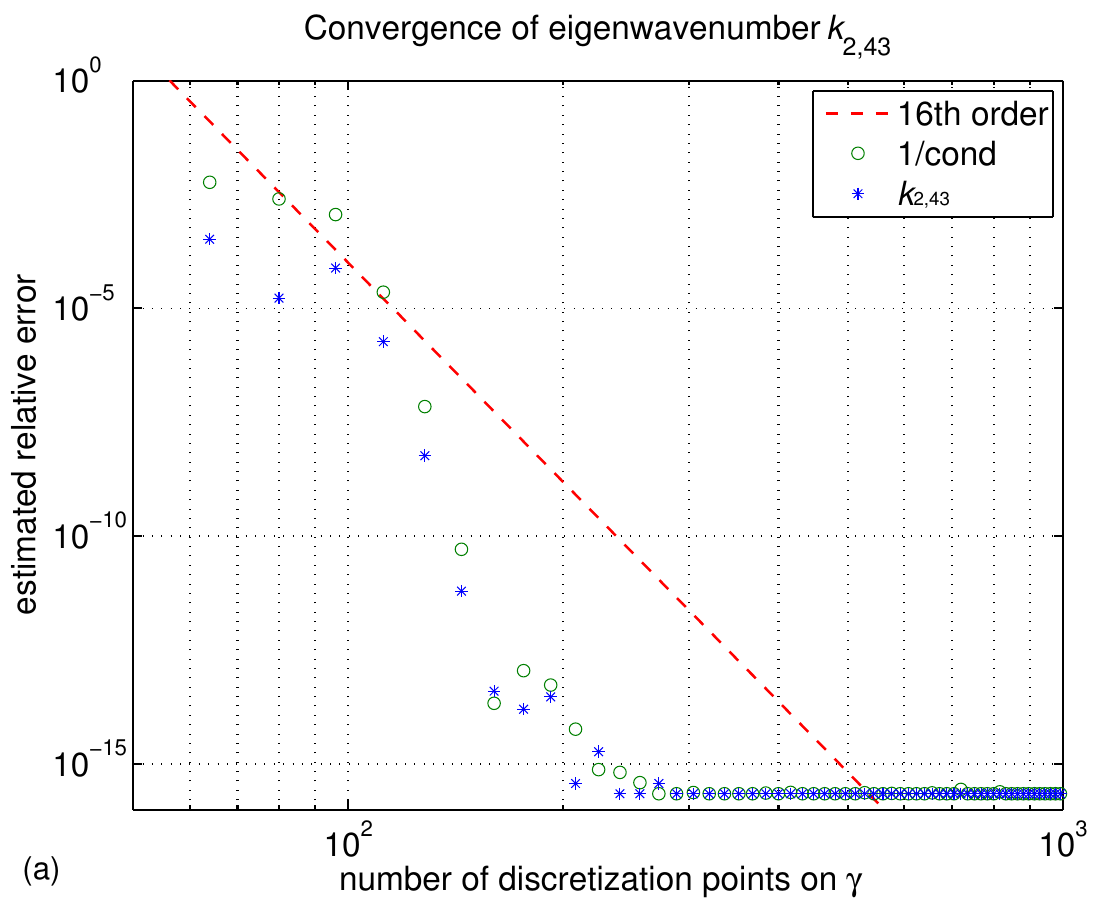}
\includegraphics[height=50mm]{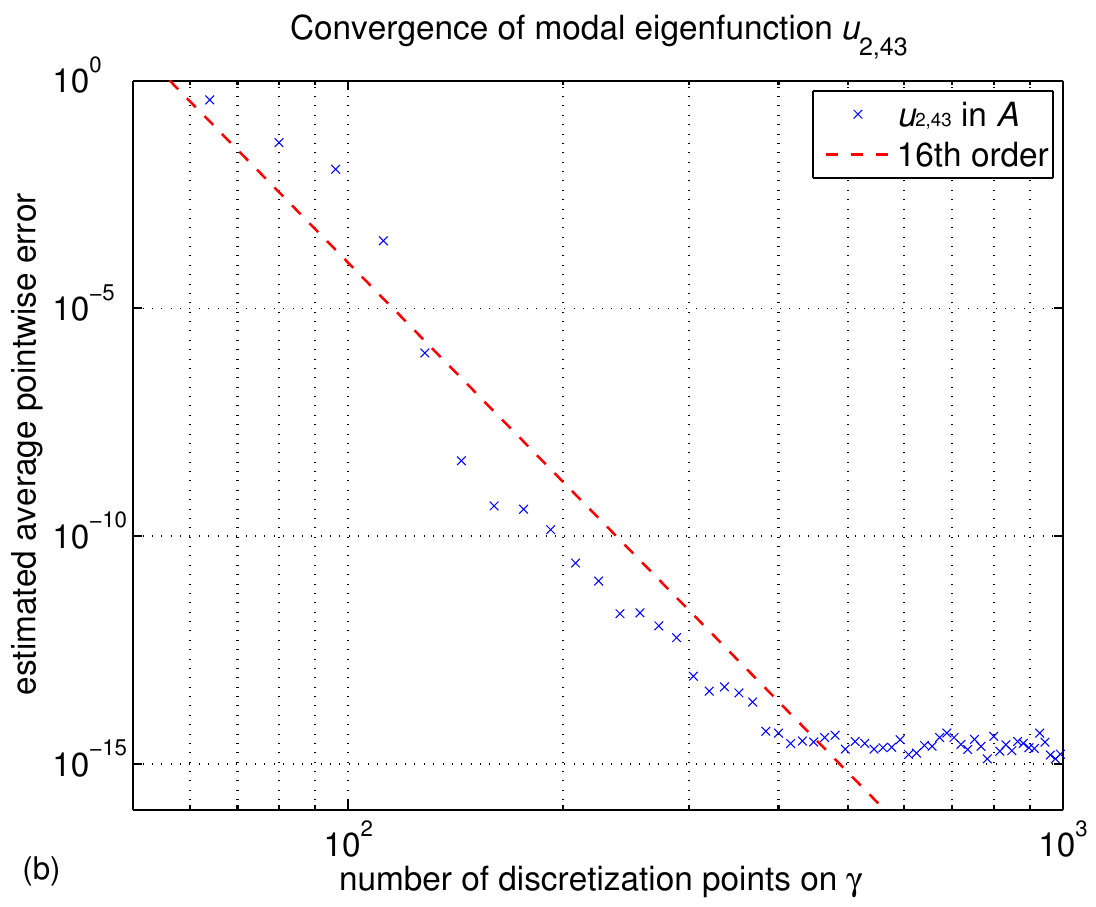}
\includegraphics[height=50mm]{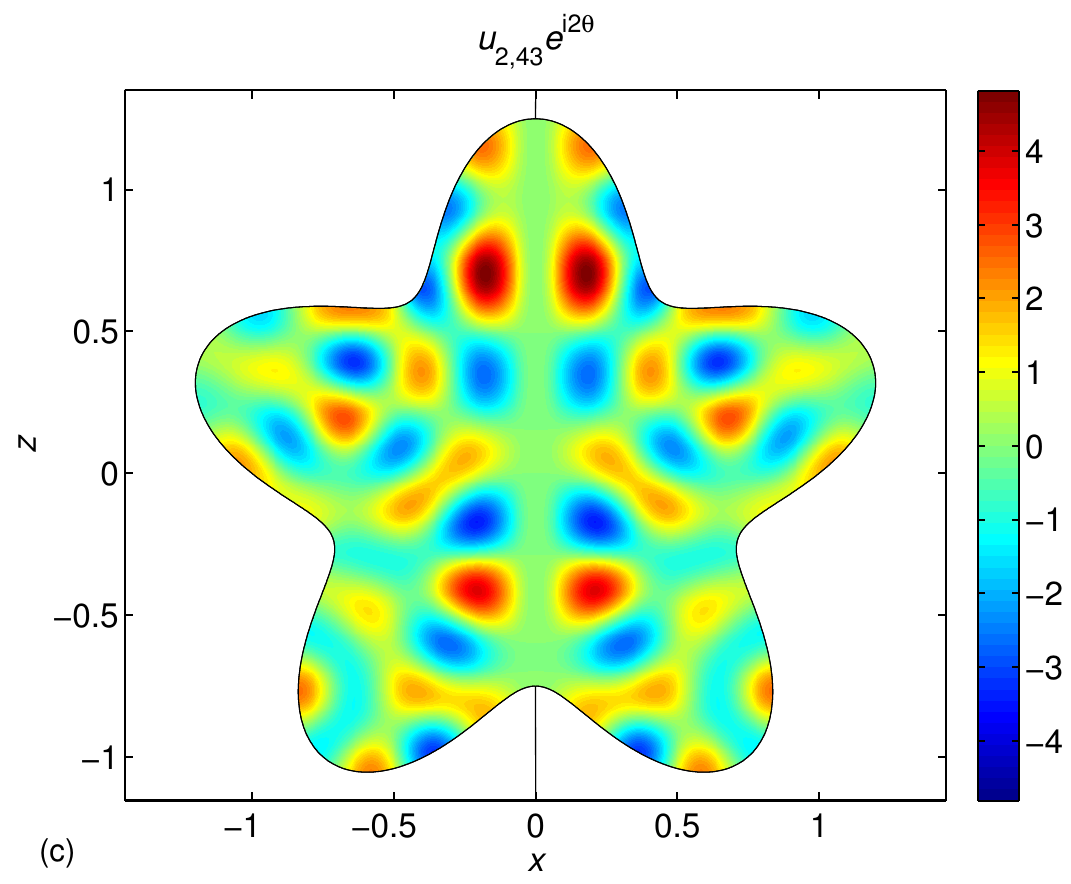}
\includegraphics[height=50mm]{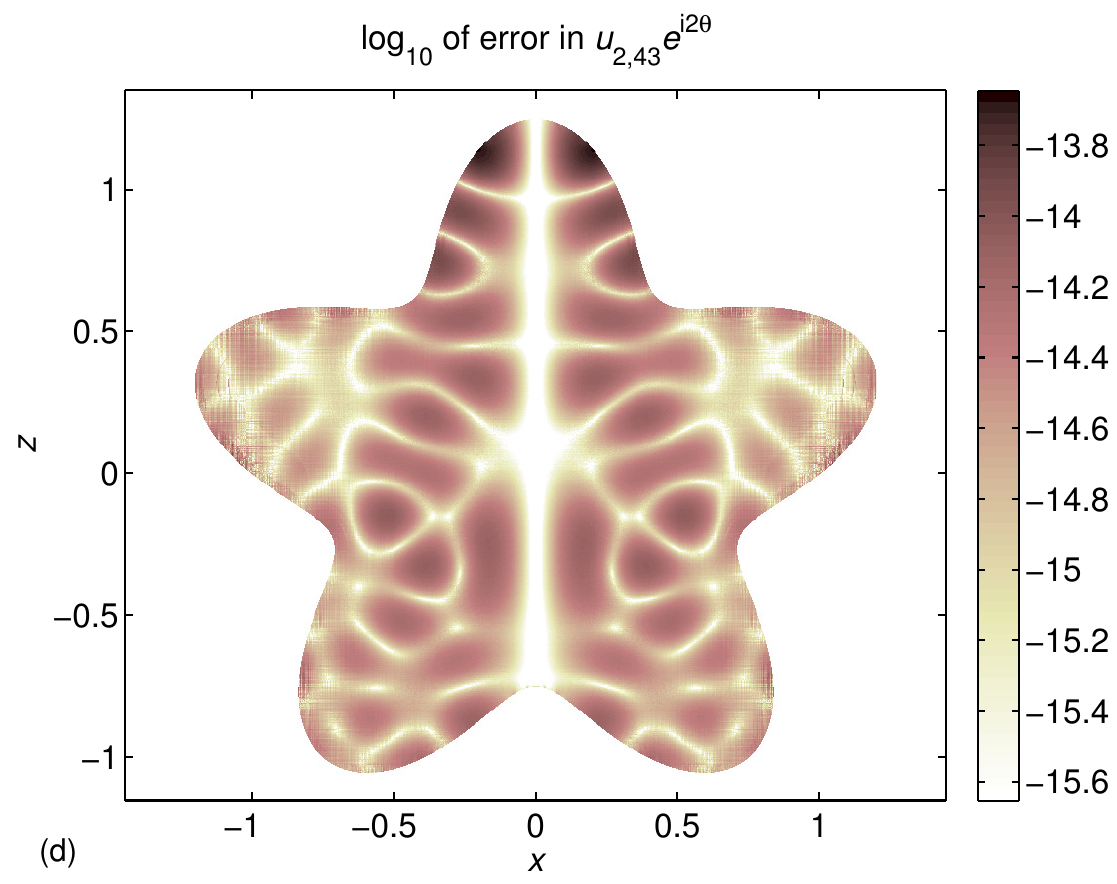}
\caption{\sf Convergence of the Neumann eigenpair $k_{2,43}$ and 
  $u_{2,43}(r)$. (a) Reciprocal condition number and error in
  $k_{2,43}$. Reference value $k_{2,43}=19.21873987061249$ is used.
  (b) Estimated average pointwise error in $u_{2,43}(r)$. (c) The
  field $u_{2,43}(r)e^{{\rm i}2\theta}$ for $\theta=0$ and
  $\theta=\pi$. (d) $\log_{10}$ of pointwise error in
  $u_{2,43}(r)e^{{\rm i}2\theta}$ for $\theta=0$ and $\theta=\pi$ with
  $560$ points on $\gamma$.}
\label{fig:n2}
\end{figure}

The results of a parallel investigation of the convergence of
$k_{2,43}$ and $u_{2,43}(r)$ are shown in Figure~\ref{fig:n2}. The
results are very similar to those for $k_{1,49}$ and $u_{1,49}(r)$ and
we conclude that computing Neumann eigenpairs is a well conditioned
problem at these wavenumbers. With our scheme it is no more difficult,
in terms of achievable accuracy, than computing simple modal fields as
in Section~\ref{sec:point}.

\section{Conclusions and outlook}

We have constructed a Fourier--Nyström discretization scheme for
second kind Fredholm integral equations with singular kernels in
axially symmetric domains and verified it numerically on modal
interior Neumann Helmholtz problems. The competitiveness of the scheme
lies in high-order convergence, high achievable accuracy, and the
ability to evaluate field solutions with uniform accuracy throughout
the computational domain. These favorable characteristics are made
possible by an explicit kernel-split, panel-based, product integration
philosophy which incorporates analytic information about integral
kernels to a higher degree than competing methods and allows for
on-the-fly computation of nearly singular quadrature rules regardless
of where target points are located relative to quadrature panels.

Another advantage of using a Fourier--Nyström scheme, over a full
three dimensional PDE eigenvalue solver for axisymmetric problems, is
that it enables easier identification and classification of
eigenfunctions. One azimuthal index $n$ is treated at a time.

One could argue that the precision offered by our scheme may not be
needed in real life applications. In acoustics, sound-hard surfaces
with homogeneous Neumann boundary conditions are only coarse models of
real surfaces. On the other hand, in electromagnetics there are
resonance problems where the mathematical models are very exact. One
example is the determination of resonance frequencies and fields in
axially symmetric superconducting cavities with ultra-high vacuum and
electro-polished surfaces. This problem is important in particle
accelerator design~\cite{Wang08} and can be modeled as a PDE
eigenvalue problem. When the electromagnetic fields are weak enough
not to affect the superconductivity, the relative error in the PDE
model can be as low as $10^{-12}$, given a generating curve $\gamma$.
Furthermore, the superconducting cavities often have corners that
subject numerical solvers to much tougher tests than the smooth
domains used in the numerical examples of the present paper. Having
already developed powerful methods for scattering problems in
non-smooth planar domains~\cite{Hels13,HelsKarl13}, we intend to
generalize our axisymmetric scheme to cope with electromagnetic
resonances and non-smooth domains in the near future. The ultimate
goal is to include our solvers in a robust particle accelerator
simulation package.

\section*{Acknowledgement}

This work was supported by the Swedish Research Council under contract
621-2011-5516.

\renewcommand{\theequation}{A.\arabic{equation}}
\setcounter{equation}{0}

\section*{Appendix A. The construction of $w_{{\rm L}j}^{\rm corr}(r_i)$}
\label{sec:AppA}

We first show how to construct the logarithmic weights $w_{{\rm
    L}j}(r_i)$ of Section~\ref{sec:Log} which occur in the
approximation of
\begin{equation}
I_p(r_i)=\int_{\gamma_p}\log|r_i-r'|G^{(1)}(r_i,r')\varrho(r')\,{\rm d}\gamma'
\label{eq:Aint1}
\end{equation}
with an expression of the form
\begin{equation}
I_p(r_i)=\sum_j G^{(1)}(r_i,r_j)\varrho_js_jw_{{\rm L}j}(r_i)\,.
\label{eq:Asum1}
\end{equation}
For simplicity we assume that $r_i$ and $r_j$ are located on the same
quadrature panel $\gamma_p$ with starting point $r(t_a)$ and end point
$r(t_b)$. Then the quadrature nodes $t_j$ and weights $w_j$,
corresponding to points $r_j$, can be expressed as
\begin{equation}
t_j=\frac{t_b+t_a}{2}+\frac{t_b-t_a}{2}\mathfrak{t}_j\,, \qquad 
w_j=\frac{t_b-t_a}{2}\mathfrak{w}_j\,,
\end{equation}
where $\mathfrak{t}_j$ and $\mathfrak{w}_j$ are $n_{\rm pt}$ nodes and
weights on the canonical panel $[-1,1]$.

Introducing $\Delta=(t_b-t_a)/2$ we rewrite~(\ref{eq:Aint1}) as
\begin{multline}
I_p(r_i)=\int_{t_a}^{t_b}G^{(1)}(r_i,r(t'))
\log\left|\frac{\Delta(r_i-r(t'))}{t_i-t'}\right|\varrho(r(t'))s(r(t'))\,{\rm d}t'\\
+\int_{t_a}^{t_b}G^{(1)}(r_i,r(t'))
\log\left|\frac{t_i-t'}{\Delta}\right|\varrho(r(t'))s(r(t'))\,{\rm d}t'\,.
\label{eq:Aint2}
\end{multline}
The first integral in~(\ref{eq:Aint2}) has a smooth integrand and is
accurately discretized as
\begin{equation}
\sum_{j\ne i}G^{(1)}(r_i,r_j)\log\left|\frac{r_i-r_j}
{\mathfrak{t}_i-\mathfrak{t}_j}\right|\varrho_js_jw_j
+G^{(1)}(r_i,r_i)\log\left|\Delta s_i\right|\varrho_is_iw_i\,.
\label{eq:Asmooth}
\end{equation}
The second integral in~(\ref{eq:Aint2}) can be transformed into
\begin{equation}
{\Delta}\int_{-1}^1G^{(1)}(r_i,r(t(\mathfrak{t}')))
\log\left|\mathfrak{t}_i-\mathfrak{t}'\right|
\varrho(r(t(\mathfrak{t}')))s(r(t(\mathfrak{t}')))\,{\rm d}\mathfrak{t}'
\end{equation}
and accurately discretized as in~\cite[Section 2.3]{Hels09}. The
result is
\begin{equation}
{\Delta}\sum_jG^{(1)}(r_i,r_j)\mathfrak{W}_{{\rm L}ij}\varrho_js_j\,,
\label{eq:Asing}
\end{equation}
where $\mathfrak{W}_{\rm L}$ is a square matrix whose entries are
$(n_{\rm pt}-1)$th degree product integration weights for the
logarithmic integral operator on the canonical panel and only depend
on the $n_{\rm pt}$ distinct nodes $\mathfrak{t}_j$.

Combining~(\ref{eq:Asmooth}) and~(\ref{eq:Asing}), the discretization
of~(\ref{eq:Aint2}) reads
\begin{multline}
I_p(r_i)=\sum_{j\ne i}G^{(1)}(r_i,r_j)\left(\log\left|r_i-r_j\right|
-\log\left|\mathfrak{t}_i-\mathfrak{t}_j\right|\right)\varrho_js_jw_j\\
+G^{(1)}(r_i,r_i)\log\left|\Delta s_i\right|\varrho_is_iw_i
+{\Delta}\sum_jG^{(1)}(r_i,r_j)\mathfrak{W}_{{\rm L}ij}\varrho_js_j\,.
\end{multline}
From~(\ref{eq:Asum1}) it is now easy to identify $w_{{\rm L}j}(r_i)$
as
\begin{equation}
w_{{\rm L}j}(r_i)=\left\{
\begin{array}{lr}
\log\left|r_i-r_j\right|w_j
-\log\left|\mathfrak{t}_i-\mathfrak{t}_j\right|w_j
+\Delta \mathfrak{W}_{{\rm L}ij}\,, & j\ne i\,,\\
\log\left|\Delta s_i\right|w_i+\Delta\mathfrak{W}_{{\rm L}ii}\,, & j=i\,.
\end{array}
\right.
\end{equation}

The definition of $w_{{\rm L}j}^{\rm corr}(r_i)$, see
Section~\ref{sec:Log}, gives
\begin{equation}
w_{{\rm L}j}^{\rm corr}(r_i)=\left\{
\begin{array}{lr}
\mathfrak{W}_{{\rm L}ij}/\mathfrak{w}_j
-\log\left|\mathfrak{t}_i-\mathfrak{t}_j\right|\,, & j\ne i\,,\\
\mathfrak{W}_{{\rm L}ii}/\mathfrak{w}_i
+\log\left|\Delta s_i\right|\,, & j=i\,.
\end{array}
\right.
\label{eq:Acorr}
\end{equation}
We observe that the weight corrections in~(\ref{eq:Acorr}) have a very
simple form, that the off-diagonal corrections do not depend on
$\gamma_p$, and that $\mathfrak{W}_{\rm L}$ only needs to be computed
and stored once. An analogous derivation for $r_i$ and $r_j$ on
neighboring panels shows that the corresponding corrections then also
depend on the relative length (in parameter) of the panels.

\renewcommand{\theequation}{B.\arabic{equation}}
\setcounter{equation}{0}
\section*{Appendix B. The construction of $w_{{\rm C}j}^{\rm cmp}(r_i)$}

The Cauchy-type singular compensation weights $w_{{\rm C}j}^{\rm
  cmp}(r_i)$ of Section~\ref{sec:LogCau} occur in the approximation of
\begin{equation}
I_p(r_i)=\int_{\gamma_p}
\frac{\mu\cdot(r_i-r')}{|r_i-r'|^2}G^{(2)}(r_i,r')\varrho(r')\,{\rm d}\gamma'
\label{eq:Bint1}
\end{equation}
with an expression of the form
\begin{equation}
I_p(r_i)=\sum_{j\ne i} G^{(2)}(r_i,r_j)\frac{\mu\cdot(r_i-r_j)}{|r_i-r_j|^2}
\varrho_js_jw_j+
\sum_j G^{(2)}(r_i,r_j)\varrho_jw_{{\rm C}j}^{\rm cmp}(r_i)\,.
\label{eq:Bsum1}
\end{equation}

Using the same notation and the same assumptions about $\gamma$,
$r_i$, and $r_j$ as in Appendix~A, we first address the construction
of product integration weights for the Cauchy operator acting on
$\varrho(r)$
\begin{equation}
J_p(\zeta_i)=\int_{\gamma_p}\frac{\varrho(\zeta')\,{\rm d}\zeta'}
{\zeta'-\zeta_i}\,.
\label{eq:BJ}
\end{equation}
Here $\zeta$ are points in the complex plane $\mathbb{C}$ which should
be identified with $r$ in $\mathbb{R}^2$ and ${\rm d}\zeta={\rm
  i}n(\zeta){\rm d}\gamma$ where the outward unit complex normal
$n(\zeta)$ corresponds to $\nu$ in $\mathbb{R}^2$. A splitting and
some change of variables give
\begin{equation}
J_p(\zeta_i)=\int_{t_a}^{t_b}\varrho(\zeta(t'))\left(
\frac{\dot{\zeta}(t')}{\zeta(t')-\zeta_i}-\frac{1}{t'-t_i}
\right)\,{\rm d}t'
+\int_{-1}^1\frac{\varrho(\zeta(t(\mathfrak{t}')))\,
{\rm d}\mathfrak{t}'}{\mathfrak{t}'-\mathfrak{t}_i}\,,
\label{eq:Bcau}
\end{equation}
where $\dot{\zeta}(t)={\rm d}\zeta(t)/{\rm d}t={\rm
  i}n(\zeta(t))s(\zeta(t))$. The first integral in~(\ref{eq:Bcau}) has
a smooth integrand and is accurately discretized as
\begin{equation}
\sum_{j\ne i}\varrho_j\left(\frac{\dot{\zeta}_jw_j}{\zeta_j-\zeta_i}
-\frac{\mathfrak{w}_j}{\mathfrak{t}_j-\mathfrak{t}_i}
\right)
+\frac{\varrho_i\ddot{\zeta}_iw_i}{2\dot{\zeta}_i}\,.
\label{eq:Bsmooth}
\end{equation}
The second integral in~(\ref{eq:Bcau}) is discretized using the
analytic method of~\cite[Section 2.1]{Hels09}, restricted to the
canonical panel. The result is
\begin{equation}
\sum_j\mathfrak{W}_{{\rm C}ij}\varrho_j\,,
\label{eq:Bsing}
\end{equation}
where $\mathfrak{W}_{\rm C}$ is a square matrix of $(n_{\rm pt}-1)$th
degree product integration weights whose entries only depend on the
$n_{\rm pt}$ distinct nodes $\mathfrak{t}_j$.
Combining~(\ref{eq:Bsmooth}) and~(\ref{eq:Bsing}), the discretization
of ~(\ref{eq:BJ}) reads
\begin{equation}
J_p(\zeta_i)=\sum_{j\ne i}\frac{{\rm i}n_j\varrho_js_jw_j}{\zeta_j-\zeta_i}
+\frac{\varrho_i\ddot{\zeta}_iw_i}{2\dot{\zeta}_i}
+\sum_{j\ne i}\varrho_j\left(\mathfrak{W}_{{\rm C}ij}
-\frac{\mathfrak{w}_j}{\mathfrak{t}_j-\mathfrak{t}_i}\right)
+\varrho_i\mathfrak{W}_{{\rm C}ii}\,,
\end{equation}
where the expression within parenthesis in the second sum can be
interpreted as a compensation weight that does not depend on
$\gamma_p$.

With access to a discretization, it remains to make a smooth
modification of the kernel in~(\ref{eq:BJ}) so that it coincides with
the kernel in~(\ref{eq:Bint1}). If, for example, $\mu=\tau$ then
kernel multiplication in~(\ref{eq:BJ}) with
\begin{equation}
-G^{(2)}(\zeta,\zeta')n(\zeta)\overline{n(\zeta')}\,,
\end{equation}
followed by taking the real part, achieves this. The compensation
weights in~(\ref{eq:Bsum1}) become
\begin{equation}
w_{{\rm C}j}^{\rm cmp}(r_i)=\left\{
\begin{array}{lr}
-(\nu_i\cdot\nu_j)(\mathfrak{W}_{{\rm C}ij}
-\mathfrak{w}_j/(\mathfrak{t}_j-\mathfrak{t}_i))\,, & \quad j\ne i\,,\\
-\mathfrak{W}_{{\rm C}ii}-\ddot{\zeta}_iw_i/(2\dot{\zeta}_i)\,, & \quad j=i\,.
\end{array}
\right.
\end{equation}
Besides a simple dependence on the unit normal, these weights share
all the desirable properties of the corrections in~(\ref{eq:Acorr}).

\begin{small}

\end{small}

\end{document}